%% file: Main.tex
\documentclass{elsarticle}
\newcommand{\english}{}
\input{Head.tex}

\input{Macros}

\title{Vapnik-Chervonenkis Dimension and Density on Johnson and Hamming Graphs }
\begin{document}
\begin{frontmatter}
\author[Computing]{Isolde Adler}
\ead{I.M.Adler@leeds.ac.uk}
\author[Maths,Computing]{Bjarki Geir Benediktsson}
\ead{B.G.Benediktsson@leeds.ac.uk}
\author[Maths]{Dugald Macpherson}
\ead{H.D.MacPherson@leeds.ac.uk}
\address[Maths]{School of Mathematics, University of Leeds}
\address[Computing]{School of Computing, University of Leeds}

\input{Abstract.tex}

\end{frontmatter}
\input{Introduction.tex}

\input{Preliminaries.tex}

\input{Johnson/Main.tex}

\input{Hamming/Main.tex}

\bibliography{Bib.bib}
\bibliographystyle{plainnat}
\input{Appendix.tex}
\if\value{TODOs}0
\else
\PackageWarning{TODO}{There are \arabic{TODOs} things to do }
\fi
\end{document}

%% file: Head.tex
\ifdefined\english
\usepackage{amsthm}
\usepackage{enumitem}

\newtheorem{definition}{Definition}[section]
\newtheorem{theorem}{Theorem}[section]
\newtheorem{lemma}[theorem]{Lemma}
\newtheorem{corollary}[theorem]{Corollary}
\else
\usepackage[icelandic]{babel}
\newtheorem{theorem}{Setning}[section]
\newtheorem{lemma}[theorem]{Hjálparsetning}

\newtheorem{corollary}[theorem]{Fylgisetning}

\newenvironment{proof}[1][Sönnun]{\begin{trivlist}
\item[\hskip \labelsep {\bfseries #1}]}{\end{trivlist}}
\newenvironment{definition}[1][Skilgreining]{\begin{trivlist}
\item[\hskip \labelsep {\bfseries #1}]}{\end{trivlist}}

\newcommand{\qed}{\nobreak \ifvmode \relax \else
      \ifdim\lastskip<1.5em \hskip-\lastskip
      \hskip1.5em plus0em minus0.5em \fi \nobreak
      \vrule height0.75em width0.5em depth0.25em\fi}

\fi

\usepackage{multicol}
\usepackage{scrextend}
\newcounter{TODOs}
\newcommand{\TODO}[1][]{\stepcounter{TODOs}\PackageWarning{mypack}{There are things to do}}

\usepackage[T1]{fontenc}
\usepackage{amsmath}
\usepackage{amsfonts}
\usepackage{sidecap}
\usepackage[utf8]{inputenc}
\usepackage{tikz}
\usepackage{listings}

\newcommand{\setR}{\ifmmode \mathbb{R} \else $\mathbb{R}$\fi}
\newcommand{\setC}{\ifmmode \mathbb{C} \else $\mathbb{C}$\fi}
\newcommand{\setN}{\ifmmode \mathbb{N} \else $\mathbb{N}$\fi}
\newcommand{\setZ}{\ifmmode \mathbb{Z} \else $\mathbb{Z}$\fi}
\newcommand{\setQ}{\ifmmode \mathbb{Q} \else $\mathbb{Q}$\fi}
\newcommand{\setE}{\ifmmode \mathbb{E} \else $\mathbb{E}$\fi}
\newcommand{\classC}{\ifmmode \mathcal{C} \else $\mathcal{C}$\fi}

\DeclareUnicodeCharacter{00AC}{\ifmmode \neg \else $\neg$\fi}
\DeclareUnicodeCharacter{025B}{\ifmmode \epsilon \else $\epsilon$}
\DeclareUnicodeCharacter{03BB}{\ifmmode \lambda \else $\lambda$\fi}
\DeclareUnicodeCharacter{03BC}{\ifmmode \mu \else $\mu$\fi}
\DeclareUnicodeCharacter{2200}{\forall}
\DeclareUnicodeCharacter{2203}{\exists}
\DeclareUnicodeCharacter{2205}{\emptyset}
\DeclareUnicodeCharacter{2208}{\in}
\DeclareUnicodeCharacter{2209}{\not\in}
\DeclareUnicodeCharacter{220B}{\ni}
\DeclareUnicodeCharacter{220C}{\not\ni}
\DeclareUnicodeCharacter{2211}{\sum}
\DeclareUnicodeCharacter{2216}{\setminus}
\DeclareUnicodeCharacter{2227}{\wedge}
\DeclareUnicodeCharacter{2228}{\vee}
\DeclareUnicodeCharacter{2229}{\cap}
\DeclareUnicodeCharacter{222A}{\cup}
\DeclareUnicodeCharacter{2260}{\neq}
\DeclareUnicodeCharacter{2264}{\leq}
\DeclareUnicodeCharacter{2265}{\geq}
\DeclareUnicodeCharacter{226E}{\not<}
\DeclareUnicodeCharacter{226F}{\not>}
\DeclareUnicodeCharacter{2270}{\not\leq}
\DeclareUnicodeCharacter{2271}{\not\geq}
\DeclareUnicodeCharacter{2282}{\subset}
\DeclareUnicodeCharacter{2283}{\supset}
\DeclareUnicodeCharacter{2284}{\not\geq}
\DeclareUnicodeCharacter{2285}{\not\geq}
\DeclareUnicodeCharacter{2286}{\subseteq}
\DeclareUnicodeCharacter{2287}{\supseteq}
\DeclareUnicodeCharacter{2288}{\not\subseteq}
\DeclareUnicodeCharacter{2289}{\not\supseteq}
\DeclareUnicodeCharacter{22C2}{\bigcap}
\DeclareUnicodeCharacter{22C3}{\bigcup}

\DeclareUnicodeCharacter{03B1}{\alpha}
\DeclareUnicodeCharacter{03B2}{\beta}
\DeclareUnicodeCharacter{03B3}{\gamma}
\DeclareUnicodeCharacter{03B4}{\delta}
\DeclareUnicodeCharacter{03B5}{\epsilon}
\DeclareUnicodeCharacter{03B6}{\zeta}
\DeclareUnicodeCharacter{03B7}{\eta}
\DeclareUnicodeCharacter{03B8}{\theta}
\DeclareUnicodeCharacter{03B9}{\iota}
\DeclareUnicodeCharacter{03BA}{\kappa}
\DeclareUnicodeCharacter{03BB}{\lambda}
\DeclareUnicodeCharacter{03BC}{\mu}
\DeclareUnicodeCharacter{03BD}{\nu}
\DeclareUnicodeCharacter{03BE}{\xi}
\DeclareUnicodeCharacter{03BF}{\omicron}
\DeclareUnicodeCharacter{03C0}{\pi}
\DeclareUnicodeCharacter{03C1}{\rho}
\DeclareUnicodeCharacter{03C3}{\sigma}
\DeclareUnicodeCharacter{03C4}{\tau}
\DeclareUnicodeCharacter{03C5}{\upsilon}
\DeclareUnicodeCharacter{03C6}{\phi}
\DeclareUnicodeCharacter{03D5}{\phi}
\DeclareUnicodeCharacter{03C7}{\chi}
\DeclareUnicodeCharacter{03C8}{\psi}
\DeclareUnicodeCharacter{03C9}{\omega}

\DeclareUnicodeCharacter{0391}{\Alpha}
\DeclareUnicodeCharacter{0392}{\Beta}
\DeclareUnicodeCharacter{0393}{\Gamma}
\DeclareUnicodeCharacter{0394}{\Delta}
\DeclareUnicodeCharacter{0395}{\Epsilon}
\DeclareUnicodeCharacter{0396}{\Zeta}
\DeclareUnicodeCharacter{0397}{\Eta}
\DeclareUnicodeCharacter{0398}{\Theta}
\DeclareUnicodeCharacter{0399}{\Iota}
\DeclareUnicodeCharacter{039A}{\Kappa}
\DeclareUnicodeCharacter{039B}{\Lambda}
\DeclareUnicodeCharacter{039C}{\Mu}
\DeclareUnicodeCharacter{039D}{\Nu}
\DeclareUnicodeCharacter{039E}{\Xi}
\DeclareUnicodeCharacter{039F}{\Omicron}
\DeclareUnicodeCharacter{03A0}{\Pi}
\DeclareUnicodeCharacter{03A1}{\Rho}
\DeclareUnicodeCharacter{03A3}{\Sigma}
\DeclareUnicodeCharacter{03A4}{\Tau}
\DeclareUnicodeCharacter{03A5}{\Upsilon}
\DeclareUnicodeCharacter{03A6}{\Phi}
\DeclareUnicodeCharacter{03A7}{\Chi}
\DeclareUnicodeCharacter{03A8}{\Psi}
\DeclareUnicodeCharacter{03A9}{\Omega}

\DeclareUnicodeCharacter{2115}{\setN}

%% file: Macros.tex
\newcommand{\margin}[1]{}
\newcommand{\isolde}[1]{\margin{\textbf{I: }#1}}

\newcommand{\VC}{\ifmmode \operatorname{VC} \else $\operatorname{VC}$\fi}
\newcommand{\vc}{\ifmmode \operatorname{vc} \else $\operatorname{vc}$\fi}
\newcommand{\symdiff}{\triangle}
\newcommand{\pow}{\mathcal P}
\newcommand{\dist}{d}
\renewcommand{\mid}{\;\colon\;}
\newcommand{\Johnson}{\ifmmode \mathcal{J}\else  $\mathcal{J}$\fi}
\newcommand{\JIS}{\ifmmode \overline{\mathcal{J}}\else $\overline{\mathcal{J}}$ \fi }
\newcommand{\Hamming}{\ifmmode \mathcal{H} \else $\mathcal{H}$\fi }
\newcommand{\HIS}{\ifmmode \overline{\mathcal{H}}\else $\overline{\mathcal{H}}$ \fi }
\usepackage{tikz} 
\usetikzlibrary{shapes}

\newif\ifConference
\Conferencefalse 
\newif\ifAppendix 
\ifConference
\Appendixfalse
\else
\Appendixtrue
\fi

%% file: Abstract.tex
\begin{abstract}
\VC-dimension and \VC-density are measures of combinatorial complexity of set systems. 
\VC-dimension was first introduced in the context of statistical learning theory, and  is tightly related to the sample complexity in PAC learning. 
\VC-density is a refinement of \VC-dimension. 
Both notions are also studied in model theory, in the context of \emph{dependent} theories. 
A set system that is definable by a formula of first-order logic with parameters has finite \VC-dimension if and only if the formula is a dependent formula. 

In this paper we study the \VC-dimension and the \VC-density of the edge relation $Exy$ on Johnson graphs and on Hamming graphs.
On a graph $G$, the set system defined by the formula $Exy$ is the vertex set of $G$ along with the collection of all \emph{open neighbourhoods} of $G$.
	We show that the edge relation has \VC-dimension at most  $4$ on Johnson graphs and at most  $3$ on Hamming graphs and these bounds are optimal.
	We furthermore show that the \VC-density of the edge relation on the class of all Johnson graphs is $2$, and on the class of all Hamming graphs the \VC-density is $2$ as well.
Moreover, we show that our bounds on the \VC-dimension carry over to the class of all induced subgraphs of Johnson graphs, and to the class of all induced subgraphs of Hamming graphs, respectively.
It also follows that the \VC-dimension of the set systems of \emph{closed neighbourhoods} in Johnson graphs and Hamming graphs is bounded.

	Johnson graphs and Hamming graphs are  well known examples of distance transitive graphs. 
	Neither of these graph classes is nowhere dense nor is there a bound on their (local) clique-width. 
	Our results contrast this by giving evidence of structural tameness of the graph classes.
\end{abstract}

%% file: Introduction.tex
\section{Introduction}

\emph{Vapnik-Chervonenkis dimension} ($\VC$-dimension) is a complexity measure of set systems. 
The related parameter \emph{$\VC$-density} provides a more refined picture of set systems that have bounded $\VC$-dimension.
First introduced in the context of statistical learning theory~\cite{VapnikCh71}, 
$\VC$-dimension also plays a key role in computational learning~\cite{Valiant84,HellersteinPRW96, GroheR17}\isolde{perhaps move grohe turan citation} as well as
in model theory~\cite{VC-I}, and it has applications in numerous areas, including graph theory~\cite{VC-dimE-P}, 
\isolde{combinatorics} computational geometry~\cite{ChazelleW89}, database theory~\cite{RiondatoACZ11}, and 
graph algorithms and complexity~\cite{Bringmann0MN16,EickmeyerGKKPRS17}.
For the definition of \VC-dimension and \VC-density see Section 2.

For fixed $k,m\in \mathbb N$ with $k\leq m$, the Johnson graph $J(m,k)$ has vertices that correspond to $k$-element subsets, of an underlying universe set of cardinality $m$, 
where two vertices are adjacent if their corresponding sets intersect in $k-1$ elements. 
Figure~\ref{fig:Johnson} shows the Johnson graph $J(4,2)$. We let $\Johnson:=\{J(m,k)\mid k,m\in \mathbb N, k\leq m\}$ denote the class of all Johnson graphs, 
and we let $\JIS$ denote the closure of $\Johnson$ under the induced subgraph relation. 
A first study of induced subgraphs of Johnson graphs has been done in~\cite{NaimiSH10}.

Hamming graphs
arise from Hamming schemes and they naturally model Hamming distance. 
For fixed $d,q\in \mathbb N$, let $S$ be a set with $|S|=q$. The Hamming graph $H(d,q)$ has vertex set $S^d$, where two vertices are adjacent if they differ in precisely one coordinate. 
Figure~\ref{fig:Hamming} shows the Hamming graph $H(3,2)$.
We let 
$\Hamming:=\{H(d,q)\mid d,q\in \mathbb N\}$ denote the class of all Hamming graphs, and we let $\HIS$
denote the closure of $\Hamming$ under induced subgraphs. The class $\HIS$ has been characterized in~\cite{KlavzarP05} via certain edge labellings.
The classes $\Johnson, \JIS,\Hamming$, and $\HIS$ admit arbitrarily large cliques as subgraphs, but nevertheless come with a highly regular structure. 

Johnson graphs and Hamming graphs are graphs of high regularity. 
They feature in different areas of computer science and mathematics, including coding theory, algebraic graph theory and model theory.
Johnson graphs also appear in László Babai's  algorithm for solving the graph isomorphism problem in quasipolynomial time~\cite{GraphIsoQuasiPoly}, where they constitute the `hard case'.

Our motivation for this work is multifaceted largely stemming from algorithmic graph theory, permutation group theory, and model theory as mentioned below.
In algorithmic graph theory structural \emph{tameness} is often linked to good algorithmic properties.
Many problems on graphs, that are algorithmically hard (e.g.\ NP-hard) in general, can be solved efficiently on classes of graphs having a \emph{tame} structure, 
such as graphs of bounded tree-width~\cite{Courcelle1990}, planar graphs, graphs excluding a fixed minor, and
nowhere dense classes of graphs~\cite{NesetrilM11a}. 
Nowhere dense classes of graphs generalise the previously mentioned classes, and in~\cite{GroheKS14} it was shown that on nowhere dense classes of graphs, every problem expressible in first-order logic is fixed-parameter tractable.
All of these classes are sparse. In particular, they cannot contain arbitrarily large cliques.
However, intuitively, cliques contain about as much information as independent sets. In~\cite{CourcelleO00}, clique-width was introduced to address this (the class of all cliques has clique-width $2$), and this was further generalised to graph classes of bounded local clique-width. 
That allowed fixed-parameter tractability for first-order logic~\cite{Grohe07}. 
Nowhere dense classes of graphs are closed under taking subgraphs, i.\,e.\ if $\mathcal C$ is a nowhere dense class of graph class, then the class obtained by closing $\mathcal C$ under subgraphs is also nowhere dense.
Graph classes of bounded (local) clique-width are closed under taking \emph{induced} subgraphs.

So-called \emph{dependent} graph classes, 
i.e.\ graph classes where every first-order formula has bounded
\VC-dimension, are a common generalisation of both nowhere dense classes of graphs~\cite{AdlerA14} and classes of bounded local clique-width~\cite{grotur04}. 
We will discuss dependent classes below and we view dependence as an interesting notion of tameness. 
The classes $\Johnson,\JIS,\Hamming,$ and $\HIS$ are somewhere dense, as arbitrarily large cliques occur as subgraphs, and they have unbounded local clique width.
Indeed, the open neighbourhood of any vertex of $J(m,k)$ induces a \emph{rook's graph} $R(m-k,k)$, cf.~Figure~\ref{fig:rook}, and the class of all rook's graphs has unbounded clique-width.
Moreover, the open $2$-neighbourhood in a Hamming graph $H(d,2)$ induces the $1$-subdivision of the complete graph on $d$ vertices, see Corollary \ref{cor:Hd2}, and it is known that the class of $1$-subdivisions of complete graphs has unbounded clique-width (cf.\ e.\,g.~\cite{AdlerBRRTV10}).
While we do not give new algorithms in this paper, our results (see Theorem~\ref{maintheorem}) provide evidence of structural tameness, despite unbounded local clique-width.

\begin{figure}
\input{Diagram.tex}
\end{figure}

Hamming graphs and Johnson graphs are regular and have large vertex transitive automorphism groups making them of particular interest in permutation group theory. 
The symmetric group $S_m$ is the full automorphism group of the Johnson graph $J(m,k)$ whenever $m≠2k$, 
and the wreath product $S_q {\rm wr} S_d$ is the full automorphism group of the Hamming graph $H(d,q)$. 
In both cases these groups act {\em distance-transitively}: if $(u,v)$ and $(u’,v’)$ are pairs of vertices with $d(u,v)=d(u’,v’)$ then there is an element $g$ in the group with $g(u)=u’$ and $g(v)=v’$.
This symmetry is exploited in some of our proofs to reduce the number of cases that need to be checked.

A major theme in recent model theory has been the study of structures which are dependent, that is, in which all formulas are dependent as described below. 
Suppose that $M$ is a first-order structure over a language $L$, and $\phi(\bar{x},\bar{y})$ is an $L$-formula with $\bar{x}=(x_i)_{i=1}^n$ and $\bar{y}=(y_i)_{i=1}^m$ (we write $|\bar{x}|=n$ and $|\bar{y}|=m$). 
For any $\bar{a}\in M^m$, put $\phi(M,\bar{a}):=\{\bar{x}\in M^n:M\models \phi(\bar{x},\bar{a})\}$. 
Then $\{\phi(M,\bar{a}):\bar{a}\in M^m\}$ is a set system in $M^n$. 
This set system has finite \VC-dimension if and only if the formula $\phi(\bar{x},\bar{y})$ is \emph{dependent}, or {\em NIP} ( does not have the{\em independence property}). 
Dependent structures include structures with {\em stable} first-order theory, such as abelian groups,  separably closed fields, and free groups, {\em o-minimal} structures (such as the real field, or even the real field equipped with the exponential function), and many Henselian valued fields such as ${\mathbb Q}_p$. 
From the viewpoint of model theory, \VC-density seems to be both a more refined invariant than \VC-dimension, and to be easier to compute. This is the viewpoint developed in the papers
\cite{VC-I} and \cite{VC-II}.
For background on dependent theories see \cite{SimonNIP}.

If $\mathcal{C}$ is a class of structures in a fixed first-order language, 
then we say the formula $\phi(\bar{x},\bar{y})$ is \emph{dependent in $\mathcal{C}$} if there is $d=d_\phi\in {\mathbb N}$ such that for every $M\in \mathcal{C}$, 
the set system $\{\phi(M,\bar{a}):\bar{a}\in M^m\}$ has \VC-dimension at most $d$, and the \emph{\VC-dimension of $\phi$ on the class} is the maximum \VC-dimension, if it exists and  $\infty$ otherwise, taken as $M$ ranges through $\mathcal{C}$. 
The class $\mathcal{C}$ is {\em dependent} if all formulas are dependent in $\mathcal{C}$. 
It is known that for fixed integer $k$, the class $\{J(m,k)\mid m\in\mathbb N\}$ is dependent, because it is first-order definable in the class of all finite sets.
%
Similarly it is also known that for a fixed integer $d$, the class $\{H(d,q)\mid q∈ℕ\}$  is dependent.
The main results of this paper give tight bounds in the case that $ϕ$ is the edge relation, i.e. $ϕ(x,y)=Exy$, for the classes where \emph{both} parameters vary.
\begin{theorem}\label{maintheorem}
	The edge relation has:
\begin{itemize}
\item	$\VC$-dimension $4$ on \Johnson, the class of all Johnson graphs. 
\item	$\VC$-dimension $3$ on \Hamming, the class of all Hamming graphs.
\item	$\VC$-density   $2$ on \Johnson, the class of all Johnson graphs.
\item	$\VC$-density   $2$ on \Hamming, the class of all Hamming graphs.
\end{itemize}
\end{theorem}

We show that the \VC-dimension of the edge relation does not increase under vertex deletion, see Lemma \ref{vdel}
and hence it follows that the $\VC$-dimension of the edge relation on $\JIS$ is $4$ and
the $\VC$-dimension of the edge relation on $\HIS$ is~$3$.

It is known that boolean combinations of dependent formulas are dependent and since equality has $\VC$-dimension at most $1$ in any model it follows that any property expressible in the language of graphs without quantifiers is dependent in $\Johnson$ and $\Hamming$. 

Using the well-known connection between $\VC$-dimension and sample complexity in the probably approximately correct (PAC) model of computational learning theory, our results imply that if $\mathcal C$ is a subset of $\JIS$ (or of $\HIS$), then every concept class definable by a quantifier-free first-order formula on $\mathcal C$ is learnable with polynomial sample complexity in the PAC model, see e.\,g.~\cite{grotur04, KearnsV94}.

The techniques we use for the proofs include identifying structural graph properties and
 symmetries that allow breaking up the problem into a feasible number of cases.

In Section~\ref{Prelim} we will cover the basic concepts and notations used throughout the paper.
Section \ref{Johnson} contains the results related to Johnson graphs and 
Section \ref{Hamming} 
contains results on Hamming graphs.

%% file: Diagram.tex
\begin{minipage}{.4\textwidth}
\begin{tikzpicture} [scale=.1, v/.style= {circle,draw,fill, scale=.2pt}, foo/.style= {draw,circle,inner sep=2pt,fill}]

    \node[v] (12) at ( 0,  5) [label={above:{\tiny$\{1,2\}$}}] {};
    \node[v] (13) at ( 5,  0) [label={[label distance = 3pt]above:{\tiny$\{1,3\}$}}] {};
    \node[v] (14) at (-5,  0) [label={[label distance = 3pt]above:{\tiny$\{1,4\}$}}] {};
    \node[v] (23) at ( 20,10) [label={above:{\tiny$\{2,3\}$}}] {};
    \node[v] (24) at (-20,10) [label={above:{\tiny$\{2,4\}$}}] {};
    \node[v] (34) at ( 0,-10) [label={below:{\tiny$\{3,4\}$}}] {};
 
    \draw (12) edge (13);
    \draw (12) edge (14);
    \draw (12) edge (23);
    \draw (12) edge (24);
    \draw (13) edge (14);
    \draw (13) edge (23);
    \draw (13) edge (34);
    \draw (14) edge (24);
    \draw (14) edge (34);
    \draw (23) edge (24);
    \draw (23) edge (34);
    \draw (24) edge (34);

\end{tikzpicture}
\caption{The Johnson graph $J(4,2)$.}
\label{fig:Johnson}
\end{minipage}
\begin{minipage}{.3\textwidth}
\begin{tikzpicture} [yscale=.5,xscale=.7, v/.style= {circle,draw,fill, scale=.2pt}, foo/.style= {draw,circle,inner sep=2pt,fill}]
    \node[v] (000) at (-1,-1) [label={[label distance = 3pt]below:{\tiny$(0,0,0)$}}] {};
    \node[v] (001) at (-1, 1) [label={[label distance = 3pt]above:{\tiny$(0,0,1)$}}] {};
    \node[v] (010) at ( 1,-1) [label={[label distance = 3pt]below:{\tiny$(0,1,0)$}}] {};
    \node[v] (011) at ( 1, 1) [label={[label distance = 3pt]above:{\tiny$(0,1,1)$}}] {};
    \node[v] (100) at (-2,-2) [label={[label distance = 3pt]below:{\tiny$(1,0,0)$}}] {};
    \node[v] (101) at (-2, 2) [label={[label distance = 3pt]above:{\tiny$(1,0,1)$}}] {};
    \node[v] (110) at ( 2,-2) [label={[label distance = 3pt]below:{\tiny$(1,1,0)$}}] {};
    \node[v] (111) at ( 2, 2) [label={[label distance = 3pt]above:{\tiny$(1,1,1)$}}] {};

    \draw (000) edge (001);
    \draw (000) edge (010);
    \draw (000) edge (100);
    \draw (001) edge (011);
    \draw (001) edge (101);
    \draw (010) edge (011);
    \draw (010) edge (110);
    \draw (100) edge (101);
    \draw (100) edge (110);
    \draw (011) edge (111);
    \draw (101) edge (111);
    \draw (110) edge (111);

\end{tikzpicture}
\caption{The Hamming graph $H(3,2)$.}
\label{fig:Hamming}
\end{minipage}
\begin{minipage}{.28\textwidth}
\newcommand{\Num}{3} 
\newcommand{\Mum}{4} 

\begin{tikzpicture} [scale=.8, v/.style= {circle,draw,fill, scale=.2pt}, foo/.style= {draw,circle,inner sep=2pt,fill}]

\foreach \j in {0,...,\Mum}{%
         \foreach \i in {0,...,\Num}{
    \node[v] (h\j;\i) at ({\j},{\i}) {};
   }
}

\foreach \k in {0,...,\Mum}{%
 \foreach \j in {0,...,\Num}{%
         \foreach \i in {\j,...,\Num}{%
    \draw (h\k;\i) edge[out=-75,in=75] (h\k;\j) ;%
   }%
 }
}
\foreach \k in {0,...,\Num}{%
 \foreach \j in {0,...,\Mum}{%
         \foreach \i in {\j,...,\Mum}{%
    \draw (h\i;\k) edge[out=165,in=15] (h\j;\k) ;%
   }%
 }
}

\end{tikzpicture}
\caption{The rook's graph $R(5.4)$.}
\label{fig:rook}
\end{minipage}

%% file: Preliminaries.tex
\section{Preliminaries}
\label{Prelim}

We let $\mathbb{N}$ denote the set of natural numbers including $0$. 
For two sets $X$ and $Y$ we use $X\symdiff Y$ to denote the symmetric difference of $X$ and $Y$ i.e. $X\symdiff Y=(X∪Y)∖(X∩Y)$. 
We use $\pow(X)$ to denote the power set of $X$. We call $|X|$ the \emph{size} of $X$.
For $k\in \mathbb N$ we let $\binom{X}{k}$ denote the set of all $k$-element subsets of $X$, i.\,e.\ $\binom{X}{k}=\{u\subseteq X\mid |u|=k\}.$

\medskip\noindent \textbf{$\VC$-dimension and $\VC$-density.}
\begin{definition}
A \emph{set system} is a pair $(X,\mathcal{S})$ consisting of a \emph{universe} set $X$ and a family $\mathcal{S}\subseteq \pow(X)$ of subsets of $X$. 
\end{definition}
Set systems are sometimes also referred to as \emph{hypergraphs} or \emph{range spaces}.
\begin{definition}
Let $(X,\mathcal{S})$ be a set system and  $A\subseteq X$ be a set.
We say that $A$ is \emph{shattered} by $\mathcal{S}$ if the class of intersections of sets in $\mathcal{S}$ with $A$ is the full powerset of $A$, i.e.\ if
$\{A\cap W\mid W\in \mathcal{S} \}=\pow(A)$.
\end{definition}
\begin{definition}
We define the \emph{shatter function} $\pi_{\mathcal{S}}: \setN \to \setN$ as 
$$
\pi_\mathcal{S}(n) := \max  \big\{ | \{S\cap A \mid S\in \mathcal{S}\} | : A\subseteq X, |A|= n\big\}.
$$
\end{definition}
We use a slight abuse of notation and say that a set $A$ is \emph{maximally shattered} for size $n$ if $|A|=n$ and $\pi_\mathcal{S}(n)=|\{S∩A\mid S∈\mathcal{S}\}|$.
\begin{definition}
The \emph{\VC-dimension} of a set system  $(X,\mathcal{S})$ is 
$$
\VC((X,\mathcal{S})) = \begin{cases}
\sup\{n \in \setN\cup\{\infty\}:  \substack{X  \text{ has a subset of  size }\\n \text{ shattered by } \mathcal{S} }\}  &\text{ if } \mathcal{S}\neq\emptyset \\
-\infty &\text{ if } \mathcal{S}=\emptyset.
\end{cases}
$$
\end{definition}
In our work we expand the above concepts to apply to classes of finite graphs in the following way.
For a class $\mathcal{C}$ of set systems the \VC-dimension of the class is
$\VC(\mathcal{C})=\sup \{\VC(X,\mathcal{S}):(X,\mathcal{S})\in \mathcal{C}\}$ if it exists and $\infty$ otherwise,
and the shatter function of $\mathcal{C}$ is  
$
\pi_{\mathcal{C}}(n) = \max\{\pi_{\mathcal{S}}(n)\ : (X,\mathcal{S})\in \mathcal{C}\}
$.

We observe that the shatter function is $2^n$ for $n$ smaller than the $\VC$-dimension of the set system but for any $n$ greater than the $\VC$-dimension it is bounded above  by a polynomial in $n$.
This is due to the Sauer-Shelah Lemma.
\begin{lemma}[Sauer-Shelah \cite{VC-I}]
\label{SauerShelah}
If $(X,\mathcal{S})$ has finite $\VC$-dimension $d$ then
$
\pi_\mathcal{S}(n)≤\sum_{i=0}^{d}\binom{n}{i}.
$
\end{lemma}

The bound on the degree of the polynomial derived from the $\VC$-dimension need not be tight. 
Since the degree of the polynomial gives a more precise measure of the combinatorial complexity of a set system, this gives rise to the following definition, 
which here we only give for classes of finite set systems.
\begin{definition}
For a class of $\mathcal{C}$ of  set systems, the \emph{\VC-density} of $\mathcal{C}$ is 
$$
\vc(\mathcal{C})=\begin{cases} 
\inf\{r\in \setR^{+}:\pi_{\mathcal{C}}(n) \in \mathcal{O}(n^r)\}, \text{ if } \VC(\mathcal{C}) < \infty \\
\infty \text{ otherwise.}
\end{cases}
$$
\end{definition}
Note that by the Sauer-Shelah Lemma $\vc(\mathcal{C})≤\VC(\mathcal{C})$.


\medskip\noindent\textbf{Graphs.} 
We consider simple, undirected graphs, i.e. graphs with no self-loops or parallel edges. 
A \emph{graph} $G$ is a pair $G =(V,E)$ where $V$ is the set of \emph{vertices} of $G$ and $E\subseteq \binom{V}{2}$ is the set of \emph{edges} of $G$. 
We also use $V(G)$ to denote the vertex set of $G$ and $E(G)$ to denote the edge set of $G$.
Two vertices $u$ and $v$ are \emph{adjacent}, if $\{u,v\}\in E$.
We denote by $N_G(v)$ the neighbourhood of $v$ in  $G$ i.e. the  set of vertices that are adjacent to $v$ in $G$ and  when $G$ is clear from the context we simply write $N(v)$.
Note that $v∉N(v)$.
A graph $H=(V',E')$ is an \emph{induced subgraph} of a graph $G(V,E)$, written $H=G[V']$ if $V'⊆V$, and $E'=E|_{V'}$, 
and we say that $V'$ induces $H$ as subgraph of $G$.
A \emph{complete} graph on $n$ vertices, denoted $K_n$, is a graph $(V,E)$ such that $|V|=n$ and $E=\binom{V}{2}$.
For a graph $G$ we say that  a set $A⊆V(G)$ is a \emph{clique} if it induces a complete graph. 
We say that $A$ is a \emph{maximal clique} if it is a clique and there is no vertex $v$ such that $A⊆N(v)$.
A \emph{path}  is a sequence $(v_i)_{i=0}^k$ of  pairwise distinct vertices such that $v_i$ is adjacent to $v_{i+1}$, and we say that  $k$ is the \emph{length} of the path.
The \emph{distance} from vertex $v$ to $u$, denoted $d(v,u)$, is the minimum length of a path from $v$ to $u$.
The \emph{$1$-subdivision} of a graph $G$ is the graph obtained from $G$ by replacing all edges of $G$ by (pairwise internally disjoint) paths of length $2$.

\begin{definition}
For $m,n∈ ℕ$, the 
\emph{rook's graph} $R(m,n)$ is the graph 
whose vertex set is  $R\times C$ where $|R|=m$ and $|C|=n$
and two distinct vertices $(i,j),(k,l)$ are adjacent if and only if $i=k$ or $j=l$.
For a fixed $i$ we call  
$\{(i,j)\mid j∈C\}$ the $i$-th row and
$\{(j,i)\mid j∈R\}$ the $i$-th column of $R(m,n)$.
\end{definition}

\begin{definition}[Johnson graphs]
Let $m,k∈ℕ$ with $m≥k$ and $X$ be a set with $|X|=m$.
	The \emph{ Johnson graph} $J(m,k)$ is the graph whose vertex set is $\binom{X}{k}$,
	where two vertices are adjacent if and only if their intersection has size $k-1$
	i.e. if their symmetric difference has size $2$.
	We call $X$ the underlying set of $J(m,k)$.
\end{definition}
We let $\Johnson:=\{J(m,k)\mid k,m\in \mathbb N, k\leq m\}$ denote the class of all Johnson graphs and \JIS its closure under taking induced subgraphs..

Examples of Johnson graphs include the octahedral graph $J(4,3)$ and the complete graph $K_n=J(n,1)$.
The following lemma is easy to verify.
\begin{lemma}[\cite{johnsondist}]
Let $u$ and $v$ be vertices in a Johnson graph.
	Then $\dist(u,v)={|u\symdiff v|}/{2}$.
\end{lemma}


\begin{definition}[Hamming graph]
Let $d,q\in \mathbb N$ and let $S$ a set with $|S|=q$.
The \emph{ Hamming graph} $H(d,q)$ is the graph whose vertices correspond to elements of $S^d$, where
two vertices are adjacent if they agree in all but one coordinate
\end{definition}
We let $\Hamming:=\{H(d,q)\mid d,q\in \mathbb N\}$ denote the class of all Hamming graphs and \HIS  its closure under taking induced subgraphs.
Note that  $H(2,n)=R(n,n)$.

\medskip\noindent\textbf{First-order logic of graphs.} 
The set of all formulas of first-order logic of graphs is defined recursively from
the \emph{atomic} formulas `$Exy$' and `$x=y$', where $x$ and $y$ are variables and `$Exy$' expresses
that $x$ and $y$ are joined by an edge, and it is closed under Boolean connectives $\neg, \wedge $ and $\vee$ and existential
quantification ($\exists$) and universal quantification ($\forall$) over vertices of the graph.
A formula is
\emph{quantifier free}, if it does not contain a quantifier. 
Since we study undirected graphs, for us $E$ is a binary relation that is symmetric and irreflexive.
We write $G\models \phi$ to say that the graph $G$ satisfies formula $\phi$.
In this paper we will focus on the atomic formula $Exy$, more precisely  we are looking at the set systems obtained by it.
The set system for $Exy$ in a graph  $G$ is $(V(G),\mathcal{S}_E)$ where \[\mathcal{S}_E:=\big\{\{x\mid G\models Exy\}\mid y∈V(G)\big\}=\{N(v)\mid v∈V(G)\}.\]
We say a set $A$ is \emph{shattered by the edge relation} in a graph $G$ if $A$ is a shattered in $(V(G),\mathcal{S}_E)$.
Moreover we will say the edge relation has any characteristic (\VC-dimension, shatter function, and \VC-density)  on a graph $G$ that the set system for the edge relation on $G$ has.
We write $\VC_E(G)$ for the $\VC$-dimension of the edge relation on a graph $G$.
\input{Induced.tex}

%% file: Induced.tex
\begin{lemma}
\label{vdel}
Let $G$ be a graph and $G':=G[V(G)∖ \{u\}]$ be a graph obtained from $G$ by deleting a single vertex $u$.
Then $\VC_E(G')≤\VC_E(G)$.

\begin{proof}

For the edge relation we have
$\mathcal{S}=\{N(v)\mid v∈ V(G)\}$.
If we delete a vertex $u$ the edge relation on the resulting subgraph $G'$ will give us the class
$\mathcal{S'}=\{N(v)\setminus\{u\} \mid v ∈ V(G)∖\{u\}\}$.
Now assume that $\VC_E(G)<\VC_E(G')$.
Then there exists a set $A⊆ V(G)∖\{u\}$ such that $|A|> \VC_{E}(G)$ and $A$ is shattered by $\mathcal{S'}$.
Since $u∉A$ we have that for all 
$S⊆ V(G)$ we get  $A∩ S = A∩(S\setminus\{u\})$,
 so
$\pow(A) = \{A∩ S\mid S∈ \mathcal{S'}\}⊆ \{A∩ S | S∈ \mathcal{S}\}$.
That means that $\mathcal{S}$ shatters $A$, in contradiction with
$|A|>\VC_{E}(G)$.
\end{proof}
\end{lemma}

%% file: Johnson/Main.tex
\section{Johnson Graphs}
\label{Johnson}
In this section we will present our results  on  the \VC-dimension and \VC-density  of the edge relation in Johnson graphs.

\input{Johnson/Lemmas/Main.tex}

\input{Johnson/Dimension/Main.tex}

\input{Johnson/Dimension/Corollary.tex}
\input{Johnson/Density/NewArg.tex}

%% file: Johnson/Lemmas/Main.tex
\input{Johnson/Lemmas/d0.tex}
\input{Johnson/Lemmas/d1.tex}
\input{Johnson/Lemmas/d2.tex}

\input{Johnson/Lemmas/Neighborhoods.tex}
\input{Johnson/Lemmas/Diag.tex}

%% file: Johnson/Lemmas/d0.tex
\ifConference
\ifAppendix
\setcounter{theorem}{\value{Jd0}}
\else
\newcounter{Jd0}
\setcounter{Jd0}{\value{theorem}}
\fi
\fi
\begin{lemma}
\label{d0}
Let $v$ be a vertex in the Johnson graph $J(m,k)$.
Then $N(v)$ induces the rook's graph $R(k,m-k)$ as a subgraph of $J(m,k)$.
\ifAppendix
\begin{proof}
Let $v$ be a vertex in the Johnson graph $J(m,k)$ and  without loss of generality assume $v=[1,k]∩ℕ$.
Every vertex in $N(v)$ has the form $(v∖\{a\})∪\{x\}$ where $a∈v$ and $x∈[k+1,m]∩ℕ$.
The mapping  $(v∖\{a\})∪\{x\}\mapsto(a,x-k)$ is a graph isomorphism $J(m,k)[N(v)]\to R(k,m-k)$.


\end{proof}
\else
\label{d0}
\fi
\end{lemma}

%% file: Johnson/Lemmas/d1.tex
\begin{lemma}
\label{d1}
Let $v$ and $w$ be vertices in a Johnson graph with
$d(v,w)=1$.
Write 
$w=(v∖\{a\})∪\{x\}$.
Then we have
$u∈N(v)∩N(w)$
if and only if
$u=(v∖\{c\})∪\{z\}$
with 
exactly one of 
$c=a$ or $z=x$.
\ifAppendix
\begin{proof}

Assume $u∈N(v)∩N(w)$. 
Then since $d(v,u)=1$ we must have $u=(v∖\{c\})∪\{z\}$ for some $c$ and $z$.
Now assume $c≠a$ and $z≠x$. 
Then we have $u\symdiff w=\{a,c,x,z\}$ so $|u\symdiff v|=4$ contradicting that $d(u,w)=1$.
So we must have either $c=a$ or $z=x$.

Conversely assume $u=(v∖\{a\})∪\{z\}$ with $x≠z$. 
Then $u∩v=v∖\{a\}$ which has size $k-1$ so $u∈N(v)$.
Also $u∩w=v∖\{a\}$ which has size $k-1$ so $u∈N(w)$.
Thus we have $u∈N(v)∩N(w)$.

Assume $u=(v∖\{c\})∪\{x\}$. 
Then $u∩v=v∖\{c\}$ which has size $k-1$ so $u∈N(v)$.
Also $u∩w=(v∖\{a,c\})∪\{x\}$ which has size $k-1$ so $u∈N(w)$.
Thus we have $u∈N(v)∩N(w)$.

Note that if we have both $c=a$ and $z=x$ then $u=w$ in contradiction with $Euw$.
\end{proof}
\else
\label{d1}
\fi
\end{lemma}

%% file: Johnson/Lemmas/d2.tex
\begin{lemma}
\label{d2}
Let $v$ and $w$ be  vertices in a Johnson graph with
$d(v,w)=2$.
We can write $w=(v∖\{a,b\})∪\{x,y\}$.
Then we have 
$u∈N(v)∩N(w)$
if and only if
$u=(v∖\{c\})∪\{z\}$
with $c∈\{a,b\}$ and $z∈\{x,y\}$.
\ifAppendix
\begin{proof}
Assume $u∈N(v)∩N(w)$. Then since $d(v,u)=1$ we must have
$u=(v∖\{c\})∪\{z\}$ for some $c∈v$ and $z∉v$.

Now assume $c∉\{a,b\}$.
Then we have 
$u\symdiff w⊇\{a,b,c\}$ contradicting that $|u\symdiff w|=2$.

Similarly
$z∈\{x,y\}$ as otherwise we have  
$u\symdiff w⊇\{x,y,z\}$ in contradiction with $|u\symdiff w|=2$.
So we must have $c∈\{a,b\}$ and $z∈\{x,y\}$.

Conversely assume $u=(v∖\{c\})∪\{z\}$ with  $c∈\{a,b\}$ and $z∈\{x,y\}$.
Assume without loss of generality $u=(v∖\{a\})∪\{x\}$. 
Then $u∩v=v∖\{a\}$ which has size $k-1$ so $u∈N(v)$.
Also
$u∩w=(v∖\{a,b\})∪\{x\}$ which has size $k-1$ so $u∈N(w)$.
Thus we have $u∈N(v)∩N(w)$.
\end{proof}
\else
\label{d2}
\fi
\end{lemma}

%% file: Johnson/Lemmas/Neighborhoods.tex
\begin{lemma}
\label{JohnsonIntersect}
Let $u$ and $v$ be vertices in the Johnson graph $J(m,k)$ then

$$
|N(u)∩N(v)|=
\begin{cases}
	k(m-k)	&\text{ if } d(u,v)=0\\
	m-1	&\text{ if } d(u,v)=1\\
	4	&\text{ if } d(u,v)=2\\
	0	&\text{ if } d(u,v)≥3
\end{cases}
$$

\begin{proof}
This follows immediately from Lemmas  \ref{d0},\ref{d1}, and \ref{d2} 
\end{proof}
\end{lemma}

%% file: Johnson/Lemmas/Diag.tex
\begin{lemma}
Let $A$ be a set of vertices in a Johnson graph shattered by the edge relation and assume $|A|≥4$. 
Then there do not exist three vertices in $A$ pairwise at distance $2$ from each other.
\label{Diag} 
\begin{proof}

Let $v$ be a vertex such that $A⊆N(v)$ and $A$ contains three vertices that are pairwise of distance $2$ from each other.
That is to say we have  
$
(v∖\{a\})∪\{x\}∈A$,$
(v∖\{b\})∪\{y\}∈A$,$
(v∖\{c\})∪\{z\}∈A
$ where $a,b,c,x,y,z$ are all distinct. 

Let $w$ be a vertex such that 
$N(w)∩A=\{
(v∖\{a\})∪\{x\},(v∖\{b\})∪\{y\},(v∖\{c\})∪\{z\}
\}$.

If $d(v,w)=1$ we can write $w=(v∖\{a_1\})∪\{x_1\}$ by Lemma \ref{d1}.
We know that 
since $(v∖\{a\})∪\{x\}∈N(w)$ we have $a_1=a$ or $x_1=x$. 

	Assume $a_1 = a$. Then since $(v∖\{b\})∪\{y\}∈N(w)$ we must have $x_1=y$, 
	so we have $w=(v∖\{a\})∪\{y\}$.
	However
	$(v∖\{c\})∪\{z\}∉N((v∖\{a\})∪\{y\})$ in contradiction to  
	$N(w)∩A=\{
(v∖\{a\})∪\{x\},(v∖\{b\})∪\{y\},(v∖\{c\})∪\{z\}
\}$.

	Alternatively  assume $x_1 = x$. Then since $(v∖\{b\})∪\{y\}∈N(w)$
	we have  $a_1=b$ so we have $w=(v∖\{b\})∪\{x\}$.
	However
	$(v∖\{c\})∪\{z\}∉N((v∖\{b\})∪\{x\})$ in contradiction to  
	$N(w)∩A=\{
(v∖\{a\})∪\{x\},(v∖\{b\})∪\{y\},(v∖\{c\})∪\{z\}
\}$.

So we must have $d(v,w)=2$ and write $w=(v∖\{a_1,a_2\})∪\{x_1,x_2\}$. By Lemma \ref{d2} we know that
since $(v∖\{a\})∪\{x\}∈N(w)$ we have $a∈\{a_1,a_2\}$ and $x∈\{x_1,x_2\}$.
 Without loss of generality we assume $a_1=a$ and $x_1=x$.

Similarly, since
 $(v∖\{b\})∪\{y\}∈N(w)$, we have   $b∈\{a,a_2\}$ and $y∈\{x,x_2\}$ 
so we have $w=(v∖\{a,b\})∪\{x,y\}$. But then $(v∖\{c\})∪\{z\}∉N(w)$, a contradiction. 

\end{proof}
\end{lemma}

%% file: Johnson/Dimension/Main.tex
\input{Johnson/Dimension/Theorem.tex}
\input{Johnson/Dimension/Tight.tex}

%% file: Johnson/Dimension/Theorem.tex
\begin{theorem}
\label{VC-dimTheorem}
The \VC-dimension of the edge relation in a Johnson graph is at most $4$.
\begin{proof}
The proof goes through a series of cases demonstrating that no vertex set of size $5$ in a Johnson graph can be shattered.
We rely on the fact that every set $A$ shattered by the edge relation must have $A⊆N(v)$ for some vertex $v$ and that
 every subset of a shattered set is also shattered which allows us to drastically reduce the number of cases we need to check.

Observe that in $J(m,k)$ we can pick an element of the underlying set and the set of all vertices not containing that element induces $J(m-1,k)$ as a subgraph of $J(m,k)$
and the set of all vertices containing that element induces $J(m-1,k-1)$. 
Thus we can assume $m$ and $k$ to be arbitrarily large and since by Lemma \ref{vdel} taking induced subgraphs can only decrease the $\VC$-dimension,  our argument then holds for all $m$ and $k$.


We will start by computing the number of configurations that can be obtained by picking $4$ vertices out of $N(v)$.
Formally the configurations, which we label Case $I$ - Case $XVI$, are the orbits of the group of automorphisms fixing $v$ in its action on $4$ element subsets of $N(v)$. 
There are $16$ and out of those $8$ are  shattered by the edge relation and $8$ are not.
We will then go through them one by one.
For those cases that are not shattered by the edge relation we will give a proof of why they are not shattered,
and in the shattered cases,
 we will  demonstrate that whichever way we choose a fifth vertex to add to those collections we will always end up with a set that is not shattered by the edge relation.

Let $A$ be a set of vertices in a Johnson graph with  $|A|=4$, and $v$ be a vertex such that $A⊆N(v)$.
\input{Johnson/Dimension/ExplainCount.tex}

We now have 16 cases and will go through them one by one demonstrating that in each case 
either $A$ can not be shattered or that adding a fifth vertex to $A$ will always result in a set that cannot be shattered.
When proving a configuration does not shatter we have to prove that there exists a subset $B⊆A$ such that there exists no $w$ for which $N(w)∩A=B$.
In all cases we will have $B≠A$ so we have to check the cases $d(v,w)=1$ and $d(v,w)=2$.

When we have to add a fifth vertex we will have to check every possible combination of $\sim_a$ and $\sim_x$ between the fifth vertex and the previous four vertices, up to a relabeling
of the $x_i$ and $a_i$.
In these subcases we will often simply observe that the fifth vertex along with $3$ of the original $4$ vertices is identical to a case which is separately proved not to shatter.


We will for each case give a diagram showing those vertices of $N(v)$ we are taking to be in $A$ arranged in rows and columns as they would be in the rook graph induced by $N(v)$.
In those cases where we do not give a proof that the four vertices selected cannot form a shattered set we will have a choice of how to pick our fifth vertex to add to $A$.
The fifth vertex we will label with the associated subcase rather than $v_5$ to avoid confusion and save space on the diagrams.
Note that row permutations just correspond to relabeling of the equivalence classes of $\sim_x$ and
column permutations correspond to relabeling of equivalence classes of $\sim_a$.

\ifConference
Here we give just $2$ of these cases, 
one where we can demonstrate the $4$ vertices cannot be shattered and 
one where we demonstrate they cannot be a subset of any $5$ vertex shattered set.
The remaining cases are proved similarly. The can be found in Appendix~\ref{app:Cases}.
\fi
\input{Johnson/Dimension/Cases/Cases.tex}

\end{proof}
\end{theorem}

%% file: Johnson/Dimension/ExplainCount.tex
Let  $v_i=(v∖\{a_i\})∪\{x_i\}$ for $i∈\{1,2,3,4\}$ be the four vertices of $A$ .
Let $\sim_x$ be the equivalence relation $v_i\sim_xv_j$ if and only if $x_i=x_j$ 
and 
 $\sim_a$ be the equivalence relation $v_i\sim_av_j$ if and only if $a_i=a_j$.
Note that if we have  $v_i\sim_xv_j$ and $v_i\sim_av_j$  then $v_i=v_j$ and by our assumption that the four vertices are distinct we have  $i=j$.

There are $5$ ways, up to permutation, to split a set of size $4$ into equivalence classes. 
These correspond to the  ways of summing up to $4$.
Not every combination of equivalence classes for $\sim_a$ and $\sim_x$ is possible.
We will now look at each of  the ways $\sim_x$ can split $A$  and give the available 
ways for $\sim_a$ to split $A$.
Note that the equivalence classes of $\sim_a$ and $\sim_x$ correspond to the columns and rows of the rook's graph induced by $N(v)$.
We now look at each of the different ways of summing up to $4$.
\begin{description}[align=left]
\item[$4$]
In this case we have  $x_1=x_2=x_3=x_4$ and we therefore must have $a_i≠a_j$ whenever $i≠j$.
This means $\sim_a$ has $4$ equivalence classes of size $1$.
This gives us Case $IX$.

\item[$3+1$]
Without loss of generality we assume $x_1=x_2=x_3≠x_4$. 
Then there are two ways for $\sim_a$ to split $A$ into equivalence classes.
It can either have $2+1+1$ or $1+1+1+1$ as the partition.
In the former case we can assume without loss of generality that 
$a_1=a_4$ and this yields Case 
$X$.
In the latter we have $a_i≠a_j$ whenever $i≠j$ and this gives us Case 
$I$.

\item[$2+2$]

Without loss of generality we assume $x_1=x_2≠x_3=x_4$.
Note that this implies $a_1≠a_2$ and $a_3≠a_4$.
We now have three ways that $\sim_a$ can split $A$ into equivalence classes.

	$\mathbf{2+2}$
	We assume without loss of generality $a_1=a_3$ and $a_2=a_4$, 
	giving us Case $II$.

	$\mathbf{2+1+1}$
	We assume without loss of generality $a_1=a_3≠a_2$, $a_1≠a_4$ and  $a_2≠a_4$.
	This gives us Case $XI$.

	$\mathbf{1+1+1+1}$
	We have $a_i≠a_j$ whenever $i≠j$, 
	yielding Case $XII$.

\item[$2+1+1$]

Without loss of generality we assume $x_1=x_2≠x_3≠x_4$ and additionally assume $x_4≠x_1$.
We can have four ways for $\sim_a$ to split $A$ into equivalence classes.

	$\mathbf{3+1}$
	Without loss of generality we can assume $a_1=a_3=a_4≠a_2$.
	This is Case 
	$XIII$.

	$\mathbf{2+2}$
	Without loss of generality we can assume 
	$a_1=a_3$ and $a_2=a_4$.
	This is Case  
	$XIV$.
	
	\TODO[Skrifa betur]
	$\mathbf{2+1+1}$	
	In this instance we have two ways of grouping the vertices with $\sim_a$ 
	that are not equivalent with relabeling.
	
	By making $a_1=a_3$ we get  Case
	$III$.
	 
	By making $a_3=a_4$ we get Case 
	$IV$.
	
	$\mathbf{1+1+1+1}$
	We have $a_i≠a_j$ whenever $i≠j$, giving us Case 
	$V$.

\item[$1+1+1+1$]

Here we have $x_1,x_2,x_3,x_4$ all distinct.
We can have four ways for $\sim_a$ to split $A$ into equivalence classes.

	$\mathbf{4}$
	Here we have  $a_1=a_2=a_3=a_4$ which is Case
	$XV$.
	
	$\mathbf{3+1}$
	Without loss of generality we may assume $a_1=a_2=a_3≠a_4$,
	giving us Case 
	$VI$.
	
	$\mathbf{2+2}$
	Without loss of generality we can assume $a_1=a_2≠a_3=a_4$ which yields Case
	$XVI$.

	$\mathbf{2+1+1}$
	Without loss of generality we can assume $a_1=a_2≠a_3≠a_4$ and $a_1≠a_4$ which gives us Case
	$VII$.
	
	$\mathbf{1+1+1+1}$
	We have $a_i≠a_j$ whenever $i≠j$. 
	This gives us Case 
	$VIII$.
\end{description}

%% file: Johnson/Dimension/Cases/Cases.tex


\begin{enumerate}[leftmargin=0pt,label=\textbf{Case \Roman*}]
%
%
\item
\hspace{0pt}
\input{Johnson/Dimension/Cases/Case3.tex}
%
%
\item 
\hspace{0pt}
\input{Johnson/Dimension/Cases/Case4.tex}

\item 
\hspace{0pt}
\input{Johnson/Dimension/Cases/Case9.tex}

\item  
\hspace{0pt}
\input{Johnson/Dimension/Cases/CaseA.tex}

%
\item
\hspace{0pt}
\input{Johnson/Dimension/Cases/CaseB.tex}

\item
\hspace{0pt}
\input{Johnson/Dimension/Cases/CaseD.tex}

\item
\hspace{0pt}
\input{Johnson/Dimension/Cases/CaseF.tex}

\item
\hspace{0pt}
\input{Johnson/Dimension/Cases/CaseG.tex}

\end{enumerate}

The remaining cases shatter, so we look at the different ways a fifth vertex can be added to the collection and demonstrate that the result cannot be a shattered set.

\begin{enumerate}[leftmargin=0pt,label=\textbf{Case \Roman*}]
\setcounter{enumi}{8}
%
%
%
\item
\hspace{0pt}
\input{Johnson/Dimension/Cases/Case1.tex}
%
%
\item
\hspace{0pt}
\input{Johnson/Dimension/Cases/Case2.tex}

%
%
\item
\hspace{0pt}
\input{Johnson/Dimension/Cases/Case5.tex}

%
%
\item
\hspace{0pt}
\input{Johnson/Dimension/Cases/Case6.tex}

%
\item
\hspace{0pt}
\input{Johnson/Dimension/Cases/Case7.tex}

%
\item
\hspace{0pt}
\input{Johnson/Dimension/Cases/Case8.tex}

\item
\hspace{0pt}
\input{Johnson/Dimension/Cases/CaseC.tex}

\item
\hspace{0pt}
\input{Johnson/Dimension/Cases/CaseE.tex}

\end{enumerate}

%% file: Johnson/Dimension/Cases/Case3.tex
%
%
\begin{minipage}[t]{.5\textwidth}
\begin{align*}
	v_1 = (v∖\{a_1\})∪\{x_1\}\\
	v_2 = (v∖\{a_2\})∪\{x_1\}\\
	v_3 = (v∖\{a_3\})∪\{x_1\}\\
	v_4 = (v∖\{a_4\})∪\{x_2\}\\
\end{align*}
\end{minipage}
\begin{minipage}[t]{.5\textwidth}
\begin{tikzpicture}[yscale=-1,thick,baseline={([yshift=9ex]current bounding box.center)}]
\node[draw,circle] (A) at ( 0, 0) {$v_1$};
\node[draw,circle] (B) at ( 1, 0) {$v_2$};
\node[draw,circle] (C) at ( 2, 0) {$v_3$};
\node[draw,circle] (D) at ( 3, 1) {$v_4$};
\end{tikzpicture}
\end{minipage}
	Let $w$ be such that
	$N(w)∩A=\{v_2,v_3,v_4\}$.
	We have 2 cases. 
	\begin{enumerate}
		\item $w=(v∖\{a\})∪\{x\}$.
		Since we have to exclude $v_1$ from $N(w)$ we must by Lemma \ref{d1} have that $a≠a_1$ and $x≠x_1$. 
		So in order to have $v_2∈N(w)$ we must  have $a=a_2$
		and in order to have $v_3∈N(w)$ we must  have $a=a_3$.
		But then $a_2=a_3$  in contradiction with $v_1≠v_2$.
		
		\item $w=(v∖\{a,b\})∪\{x,y\}$. From Lemma \ref{d2} we get that
	  		$v_2∈N(w)$ yields $a_2∈\{a,b\}$ and $x_1∈\{x,y\}$;
			$v_3∈N(w)$ yields $a_3∈\{a,b\}$ and $x_1∈\{x,y\}$;
			$v_4∈N(w)$ yields $a_4∈\{a,b\}$ and $x_2∈\{x,y\}$.
		Thus $\{a_2,a_3,a_4\}⊆\{a,b\}$, contradicting that $a_2,a_3,a_4$ are all distinct.

	\end{enumerate}

%% file: Johnson/Dimension/Cases/Case4.tex
%
%
\begin{minipage}[t]{.5\textwidth}
\begin{align*}
	v_1 = (v∖\{a_1\})∪\{x_1\}\\
	v_2 = (v∖\{a_2\})∪\{x_1\}\\
	v_3 = (v∖\{a_1\})∪\{x_2\}\\
	v_4 = (v∖\{a_2\})∪\{x_2\}\\
\end{align*}
\end{minipage}
\begin{minipage}[t]{.5\textwidth}
\begin{tikzpicture}[yscale=-1,thick,baseline={([yshift=9ex]current bounding box.center)}]
\node[draw,circle] (A) at ( 0, 0) {$v_1$};
\node[draw,circle] (B) at ( 1, 0) {$v_2$};
\node[draw,circle] (C) at ( 0, 1) {$v_3$};
\node[draw,circle] (D) at ( 1, 1) {$v_4$};
\end{tikzpicture}
\end{minipage}
	Let $w$ be such that  
	$N(w)∩A=\{v_2,v_3,v_4\}$.
	We have 2 cases.
	\begin{enumerate}
		\item$w=(v∖\{a\})∪\{x\}$.
		Since $v_4∈N(w)$ we have $w≠v_1$.
		Since we have to exclude $v_1$ from $N(w)$, by Lemma \ref{d1} we must have that $a≠a_1$ and $x≠x_1$. 
		So in order to have $v_2∈N(w)$ we must  have $a=a_2$
		and in order to have $v_3∈N(w)$ we must  have $x=x_2$.
		But then $w=v_4$ in contradiction with $v_4∈N(w)$.

		\item $w=(v∖\{a,b\})∪\{x,y\}$. From Lemma \ref{d2} we get that
	  		$v_3∈N(w)$ yields $a_1∈\{a,b\}$ and $x_2∈\{x,y\}$;
			$v_2∈N(w)$ yields $a_2∈\{a,b\}$ and $x_1∈\{x,y\}$;
		hence $w=(v∖\{a_1,a_2\})∪\{x_1,x_2\}$, contradicting that $v_1∉N(w)$.
	\end{enumerate}

%% file: Johnson/Dimension/Cases/Case9.tex
\begin{minipage}[t]{.5\textwidth}
\begin{align*}
	v_1 = (v∖\{a_1\})∪\{x_1\}\\
	v_2 = (v∖\{a_2\})∪\{x_1\}\\
	v_3 = (v∖\{a_1\})∪\{x_2\}\\
	v_4 = (v∖\{a_3\})∪\{x_3\}\\
\end{align*}
\end{minipage}
\begin{minipage}[t]{.5\textwidth}
\begin{tikzpicture}[yscale=-1,thick,baseline={([yshift=9ex]current bounding box.center)}]
\node[draw,circle] (A) at ( 0, 0) {$v_1$};
\node[draw,circle] (B) at ( 1, 0) {$v_2$};
\node[draw,circle] (C) at ( 0, 1) {$v_3$};
\node[draw,circle] (D) at ( 2, 2) {$v_4$};
\end{tikzpicture}
\end{minipage}
The vertices $v_2,v_3,v_4$ are at distance $2$ from each other so by Lemma \ref{Diag} $A$ is not shattered.
%
%
%

%% file: Johnson/Dimension/Cases/CaseA.tex
\begin{minipage}[t]{.5\textwidth}
\begin{align*}
	v_1 = (v∖\{a_1\})∪\{x_1\}\\
	v_2 = (v∖\{a_2\})∪\{x_1\}\\
	v_3 = (v∖\{a_3\})∪\{x_2\}\\
	v_4 = (v∖\{a_3\})∪\{x_3\}\\
\end{align*}
\end{minipage}
\begin{minipage}[t]{.5\textwidth}
\begin{tikzpicture}[yscale=-1,thick,baseline={([yshift=9ex]current bounding box.center)}]
\node[draw,circle] (A) at ( 0, 0) {$v_1$};
\node[draw,circle] (B) at ( 1, 0) {$v_2$};
\node[draw,circle] (C) at ( 2, 1) {$v_3$};
\node[draw,circle] (D) at ( 2, 2) {$v_4$};
\end{tikzpicture}
\end{minipage}

	Let $w$ be such that  
	 $N(w)∩A=\{v_2,v_3,v_4\}$.
	We have 2 cases.
\begin{enumerate} 
		\item $w=(v∖\{a\})∪\{x\}$. 
		Since we have to exclude $v_1$ from $N(w)$ by Lemma \ref{d1} we must have that $a≠a_1$ and $x≠x_1$. 
		So in order to have $v_2∈N(w)$ we must  have $a=a_2$
		and in order to have $v_3∈N(w)$ we must  have $x=x_2$.
		But then $w=(v∖\{a_2\})∪\{x_2\}$  in contradiction with $v_4∈N(w)$.
		
		\item $w=(v∖\{a,b\})∪\{x,y\}$. From Lemma \ref{d2} we get that
	  		$v_2∈N(w)$ yields $a_2∈\{a,b\}$ and $x_1∈\{x,y\}$;
			$v_3∈N(w)$ yields $a_3∈\{a,b\}$ and $x_2∈\{x,y\}$.
		Thus $w=(v∖\{a_2,a_3\})∪\{x_1,x_2\}$ in contradiction with $v_4∈N(w)$.
\end{enumerate}

%% file: Johnson/Dimension/Cases/CaseB.tex
%
\begin{minipage}[t]{.5\textwidth}
\begin{align*}
	v_1 = (v∖\{a_1\})∪\{x_1\}\\
	v_2 = (v∖\{a_2\})∪\{x_1\}\\
	v_3 = (v∖\{a_3\})∪\{x_2\}\\
	v_4 = (v∖\{a_4\})∪\{x_3\}\\
\end{align*}
\end{minipage}
\begin{minipage}[t]{.5\textwidth}
\begin{tikzpicture}[yscale=-1,thick,baseline={([yshift=9ex]current bounding box.center)}]
\node[draw,circle] (A) at ( 0, 0) {$v_1$};
\node[draw,circle] (B) at ( 1, 0) {$v_2$};
\node[draw,circle] (C) at ( 2, 1) {$v_3$};
\node[draw,circle] (D) at ( 3, 2) {$v_4$};
\end{tikzpicture}
\end{minipage}
The vertices $v_2,v_3,v_4$ are at distance $2$ from each other so by Lemma \ref{Diag} $A$ is not shattered.
%
%

%% file: Johnson/Dimension/Cases/CaseD.tex
\begin{minipage}[t]{.5\textwidth}
\begin{align*}
	v_1 = (v∖\{a_1\})∪\{x_1\}\\
	v_2 = (v∖\{a_1\})∪\{x_2\}\\
	v_3 = (v∖\{a_1\})∪\{x_3\}\\
	v_4 = (v∖\{a_2\})∪\{x_4\}\\
\end{align*}
\end{minipage}
\begin{minipage}[t]{.5\textwidth}
\begin{tikzpicture}[yscale=-1,thick,baseline={([yshift=9ex]current bounding box.center)}]
\node[draw,circle] (A) at (0 , 0) {$v_1$};
\node[draw,circle] (B) at (0 , 1) {$v_2$};
\node[draw,circle] (C) at (0 , 2) {$v_3$};
\node[draw,circle] (D) at (1 , 3) {$v_4$};
\end{tikzpicture}
\end{minipage}
	Let $w$ be such that  
	 $N(w)∩A=\{v_2,v_3,v_4\}$.
	
	We have 2 cases.
		\begin{enumerate} 
		\item$w=(v∖\{a\})∪\{x\}$.
		Since we have to exclude $v_1$ from $N(w)$ by Lemma \ref{d1} we must have that $a≠a_1$ and $x≠x_1$ 
		so in order to have $v_2∈N(w)$ we must  have $x=x_2$
		and in order to have $v_3∈N(w)$ we must  have $x=x_3$.
		But then $x_2=x_3$,  in contradiction with $v_1≠v_2$.

		\item$w=(v∖\{a,b\})∪\{x,y\}$. From Lemma \ref{d2} we get that 
	  	$v_2∈N(w)$ yields $a_1∈\{a,b\}$ and $x_2∈\{x,y\}$;
		$v_3∈N(w)$ yields $a_1∈\{a,b\}$ and $x_3∈\{x,y\}$;
		$v_4∈N(w)$ yields $a_2∈\{a,b\}$ and $x_4∈\{x,y\}$.
		Thus we have  $\{x_2,x_3,x_4\}⊆\{x,y\}$, in contradiction with $x_2,x_3,x_4$ all being distinct.
\end{enumerate}

%% file: Johnson/Dimension/Cases/CaseF.tex
\begin{minipage}[t]{.5\textwidth}
\begin{align*}
	v_1 = (v∖\{a_1\})∪\{x_1\}\\
	v_2 = (v∖\{a_1\})∪\{x_2\}\\
	v_3 = (v∖\{a_2\})∪\{x_3\}\\
	v_4 = (v∖\{a_3\})∪\{x_4\}\\
\end{align*}
\end{minipage}
\begin{minipage}[t]{.5\textwidth}
\begin{tikzpicture}[yscale=-1,thick,baseline={([yshift=9ex]current bounding box.center)}]
\node[draw,circle] (A) at (0 , 0) {$v_1$};
\node[draw,circle] (B) at (0 , 1) {$v_2$};
\node[draw,circle] (C) at (1 , 2) {$v_3$};
\node[draw,circle] (D) at (2 , 3) {$v_4$};
\end{tikzpicture}
\end{minipage}
The vertices $v_2,v_3,v_4$ are at distance $2$ from each other so by Lemma \ref{Diag} $A$ is not shattered.
%
%
%

%% file: Johnson/Dimension/Cases/CaseG.tex
\begin{minipage}[t]{.5\textwidth}
\begin{align*}
	v_1 = (v∖\{a_1\})∪\{x_1\}\\
	v_2 = (v∖\{a_2\})∪\{x_2\}\\
	v_3 = (v∖\{a_3\})∪\{x_3\}\\
	v_4 = (v∖\{a_4\})∪\{x_4\}\\
\end{align*}
\end{minipage}
\begin{minipage}[t]{.5\textwidth}
\begin{tikzpicture}[yscale=-1,thick,baseline={([yshift=9ex]current bounding box.center)}]
\node[draw,circle] (A) at (0 , 0) {$v_1$};
\node[draw,circle] (B) at (1 , 1) {$v_2$};
\node[draw,circle] (C) at (2 , 2) {$v_3$};
\node[draw,circle] (D) at (3 , 3) {$v_4$};
\end{tikzpicture}
\end{minipage}
The vertices $v_2,v_3,v_4$ are at distance $2$ from each other so by Lemma \ref{Diag} $A$ is not shattered.
%
%

%% file: Johnson/Dimension/Cases/Case1.tex
%
%
%
\begin{minipage}[t]{.5\textwidth}
\begin{align*}
	v_1 = (v∖\{a_1\})∪\{x_1\}\\
	v_2 = (v∖\{a_2\})∪\{x_1\}\\
	v_3 = (v∖\{a_3\})∪\{x_1\}\\
	v_4 = (v∖\{a_4\})∪\{x_1\}\\
\end{align*}
\end{minipage}
\begin{minipage}[t]{.5\textwidth}
\begin{tikzpicture}[yscale=-1,thick,baseline={([yshift=9ex]current bounding box.center)}]
\node[draw,circle] (A) at ( 0, 0) {$v_1$};
\node[draw,circle] (B) at ( 1, 0) {$v_2$};
\node[draw,circle] (C) at ( 2, 0) {$v_3$};
\node[draw,circle] (D) at ( 3, 0) {$v_4$};
\node[red,draw,circle]		(D) at ( 4, 0) {a};
\node[green,draw,circle]	(D) at ( 0, 1) {b};
\node[blue,draw,circle]		(D) at ( 4, 1) {c};
\end{tikzpicture}
\end{minipage}
	This case shatters so we  take a closer look  at what configurations are obtainable by adding a fifth vertex.
	\begin{enumerate}
	\item[\textcolor{red}{a}]
	$v_5 = (v∖\{a_5\})∪\{x_1\}$.
	Let $w$ be such that $N(w)∩A=\{v_1,v_2,v_3\}$.
	Observe that $w≠v_4$ since $v_5∈N(v_4)$ so we will need an alternative $w$.
	We have $2$ cases: either $d(v,w)=1$
	 or $d(v,w)=2$.

	Let $w=(v∖\{a\})∪\{x\}$. Since we have to exclude $v_5$ from $N(w)$ then by Lemma \ref{d1} we cannot have $x=x_1$.
	So in order to have $v_1∈N(w)$ we must have $a=a_1$ but then in order to have $v_2∈N(w)$ we must have $x=x_1$, a contradiction.

	Let $w=(v∖\{a,b\})∪\{x,y\}$. 
	In order to have $v_1∈N(w),v_2∈N(w)$ and $v_3∈N(w)$, Lemma \ref{d2} gives us
	$\{a_1,a_2,a_3\}⊆\{a,b\}$, a contradiction. 	

	\item[\textcolor{green}{b}]
	$v_5 = (v∖\{a_1\})∪\{x_2\}$.
	Here $v_2,v_3,v_4,v_5$ form case $I$.
	\item[\textcolor{blue}{c}]
	$v_5 = (v∖\{a_5\})∪\{x_2\}$.
	Here $v_1,v_2,v_3,v_5$ form case $I$. 
	\end{enumerate}

%% file: Johnson/Dimension/Cases/Case2.tex
%
%
\begin{minipage}[t]{.5\textwidth}
\begin{align*}
	v_1 = (v∖\{a_1\})∪\{x_1\}\\
	v_2 = (v∖\{a_2\})∪\{x_1\}\\
	v_3 = (v∖\{a_3\})∪\{x_1\}\\
	v_4 = (v∖\{a_1\})∪\{x_2\}\\
\end{align*}
\end{minipage}
\begin{minipage}[t]{.5\textwidth}
\begin{tikzpicture}[yscale=-1,thick,baseline={([yshift=9ex]current bounding box.center)}]
\node[draw,circle] (A) at ( 0, 0) {$v_1$};
\node[draw,circle] (B) at ( 1, 0) {$v_2$};
\node[draw,circle] (C) at ( 2, 0) {$v_3$};
\node[draw,circle] (D) at ( 0, 1) {$v_4$};
\node[draw,circle,red] 		(D) at (3,0) {a};
\node[draw,circle,green]	(D) at (1,1) {b};
\node[draw,circle,blue]		(D) at (3,1) {c};
\node[draw,circle,brown]	(D) at (0,2) {d};
\node[draw,circle,gray] 	(D) at (1,2) {e};
\node[draw,circle,cyan]		(D) at (3,2) {f};
\end{tikzpicture}
\end{minipage}
	\begin{enumerate}
	\item[\textcolor{red}{a}]
	$v_5 = (v∖\{a_4\})∪\{x_1\}$. 
	Then $v_2,v_3,v_4,v_5$ forms case $I$. 
	\item[\textcolor{green}{b}]
	$v_5 = (v∖\{a_2\})∪\{x_2\}$.
	 Then $v_1,v_2,v_4,v_5$ forms case $II$. 
	\item[\textcolor{blue}{c}]
	$v_5 = (v∖\{a_4\})∪\{x_2\}$.
	 Then $v_1,v_2,v_3,v_5$ forms case $I$.
	\item[\textcolor{brown}{d}]
	$v_5 = (v∖\{a_1\})∪\{x_3\}$.
	 Then $v_2,v_3,v_4,v_5$ forms case $IV$. 
	\item[\textcolor{gray}{e}]
	$v_5 = (v∖\{a_2\})∪\{x_3\}$.
	 Then $v_1,v_3,v_4,v_5$ forms case $III$.
	\item[\textcolor{cyan}{f}]
	$v_5 = (v∖\{a_4\})∪\{x_3\}$.
	 Then $v_1,v_2,v_4,v_5$ forms case $III$.
	\end{enumerate}

%% file: Johnson/Dimension/Cases/Case5.tex
%
%
\begin{minipage}[t]{.5\textwidth}
\begin{align*}
	v_1 = (v∖\{a_1\})∪\{x_1\}\\
	v_2 = (v∖\{a_2\})∪\{x_1\}\\
	v_3 = (v∖\{a_1\})∪\{x_2\}\\
	v_4 = (v∖\{a_3\})∪\{x_2\}\\
\end{align*}
\end{minipage}
\begin{minipage}[t]{.5\textwidth}
\begin{tikzpicture}[yscale=-1,thick,baseline={([yshift=9ex]current bounding box.center)}]
\node[draw,circle] (A) at ( 0, 0) {$v_1$};
\node[draw,circle] (B) at ( 1, 0) {$v_2$};
\node[draw,circle] (C) at ( 0, 1) {$v_3$};
\node[draw,circle] (D) at ( 2, 1) {$v_4$};
\node[draw,circle,red]	(D) at ( 2, 0) {a};
\node[draw,circle,green](D) at ( 3, 0) {b};
\node[draw,circle,blue]	(D) at ( 0, 2) {c};
\node[draw,circle,brown](D) at ( 1, 2) {d};
\node[draw,circle,gray]	(D) at ( 3, 2) {e};
\end{tikzpicture}
\end{minipage}
	\begin{enumerate}
	\item[\textcolor{red}{a}]
	$v_5 = (v∖\{a_3\})∪\{x_1\}$.
	Here $v_1,v_3,v_4,v_5$ form case $II$.
	\item[\textcolor{green}{b}]
	$v_5 = (v∖\{a_4\})∪\{x_1\}$.
	Here $v_1,v_2,v_4,v_5$ form case $I$.
	\item[\textcolor{blue}{c}]
	$v_5 = (v∖\{a_1\})∪\{x_3\}$.
	Here $v_2,v_4,v_5$ all have distance $2$ from each other and thus by Lemma \ref{Diag} $A$ is not shattered.
	\item[\textcolor{brown}{d}]
	$v_5 = (v∖\{a_2\})∪\{x_3\}$.
	Here $v_1,v_4,v_5$ all have distance $2$ from each other and thus by Lemma \ref{Diag} $A$ is not shattered.
	\item[\textcolor{gray}{e}]
	$v_5 = (v∖\{a_4\})∪\{x_3\}$
	In this case $v_1,v_4,v_5$ all have distance $2$ from each other and thus by Lemma \ref{Diag} $A$ is not shattered.
	\end{enumerate}

%% file: Johnson/Dimension/Cases/Case6.tex
%
%
\begin{minipage}[t]{.5\textwidth}
\begin{align*}
	v_1 = (v∖\{a_1\})∪\{x_1\}\\
	v_2 = (v∖\{a_2\})∪\{x_1\}\\
	v_3 = (v∖\{a_3\})∪\{x_2\}\\
	v_4 = (v∖\{a_4\})∪\{x_2\}\\
\end{align*}
\end{minipage}
\begin{minipage}[t]{.5\textwidth}
\begin{tikzpicture}[yscale=-1,thick,baseline={([yshift=9ex]current bounding box.center)}]
\node[draw,circle] (A) at ( 0, 0) {$v_1$};
\node[draw,circle] (B) at ( 1, 0) {$v_2$};
\node[draw,circle] (C) at ( 2, 1) {$v_3$};
\node[draw,circle] (D) at ( 3, 1) {$v_4$};
\node[draw,circle,red]		(D) at (2,0) {a};
\node[draw,circle,green]	(D) at (4,0) {b};
\node[draw,circle,blue]		(D) at (0,2) {c};
\node[draw,circle,brown] 	(D) at (4,2) {d};
\end{tikzpicture}
\end{minipage}
	\begin{enumerate}
	\item[\textcolor{red}{a}]
	$v_5 = (v∖\{a_3\})∪\{x_1\}$.
	Then $v_1,v_2,v_4,v_5$ form case $I$.
	\item[\textcolor{green}{b}]
	$v_5 = (v∖\{a_5\})∪\{x_1\}$.
	Then $v_1,v_2,v_3,v_5$ form case $I$.
	\item[\textcolor{blue}{c}]
	$v_5 = (v∖\{a_1\})∪\{x_3\}$.
	In this case $v_1,v_3,v_4,v_5$ form case $IV$.
	\item[\textcolor{brown}{d}]
	$v_5 = (v∖\{a_5\})∪\{x_3\}$.
	In this case $v_1,v_3,v_5$ all have distance $2$ from each other and thus by Lemma \ref{Diag} $A$ is not shattered.
	\end{enumerate}

%% file: Johnson/Dimension/Cases/Case7.tex
%
\begin{minipage}[t]{.5\textwidth}
\begin{align*}
	v_1 = (v∖\{a_1\})∪\{x_1\}\\
	v_2 = (v∖\{a_2\})∪\{x_1\}\\
	v_3 = (v∖\{a_1\})∪\{x_2\}\\
	v_4 = (v∖\{a_1\})∪\{x_3\}\\
\end{align*}
\end{minipage}
\begin{minipage}[t]{.5\textwidth}
\begin{tikzpicture}[yscale=-1,thick,baseline={([yshift=9ex]current bounding box.center)}]
\node[draw,circle] (A) at ( 0, 0) {$v_1$};
\node[draw,circle] (B) at ( 1, 0) {$v_2$};
\node[draw,circle] (C) at ( 0, 1) {$v_3$};
\node[draw,circle] (D) at ( 0, 2) {$v_4$};
\node[draw,circle,red]		(D) at (2,0) {a};
\node[draw,circle,green]	(D) at (1,1) {b};
\node[draw,circle,blue]		(D) at (2,1) {c};
\node[draw,circle,brown]	(D) at (0,3) {d};
\node[draw,circle,gray]		(D) at (1,3) {e};
\node[draw,circle,cyan]		(D) at (2,3) {f};
\end{tikzpicture}
\end{minipage}
	\begin{enumerate}
	\item[\textcolor{red}{a}]
	$v_5 = (v∖\{a_3\})∪\{x_1\}$.
	Then $v_2,v_3,v_4,v_5$ form case $IV$.
	\item[\textcolor{green}{b}]
	$v_5 = (v∖\{a_2\})∪\{x_2\}$.
	Then $v_1,v_2,v_3,v_5$ form case $II$.
	\item[\textcolor{blue}{c}]
	$v_5 = (v∖\{a_3\})∪\{x_2\}$.
	In this case $v_2,v_4,v_5$ all have distance $2$ from each other and thus by Lemma \ref{Diag} $A$ is not shattered.
	\item[\textcolor{brown}{d}]
	$v_5 = (v∖\{a_1\})∪\{x_4\}$.
	Here $v_2,v_3,v_4,v_5$ form case $VI$.
	\item[\textcolor{gray}{e}]
	$v_5 = (v∖\{a_2\})∪\{x_4\}$.
	In this case $v_1,v_3,v_4,v_5$ form case $VI$.
	\item[\textcolor{cyan}{f}]
	$v_5 = (v∖\{a_3\})∪\{x_4\}$.
	In this case $v_2,v_4,v_5$ all have distance $2$ from each other and thus by Lemma \ref{Diag} $A$ is not shattered.
	\end{enumerate}

%% file: Johnson/Dimension/Cases/Case8.tex
%
\begin{minipage}[t]{.5\textwidth}
\begin{align*}
	v_1 = (v∖\{a_1\})∪\{x_1\}\\
	v_2 = (v∖\{a_2\})∪\{x_1\}\\
	v_3 = (v∖\{a_1\})∪\{x_2\}\\
	v_4 = (v∖\{a_2\})∪\{x_3\}\\
\end{align*}
\end{minipage}
\begin{minipage}[t]{.5\textwidth}
\begin{tikzpicture}[yscale=-1,thick,baseline={([yshift=9ex]current bounding box.center)}]
\node[draw,circle] (A) at ( 0, 0) {$v_1$};
\node[draw,circle] (B) at ( 1, 0) {$v_2$};
\node[draw,circle] (C) at ( 0, 1) {$v_3$};
\node[draw,circle] (D) at ( 1, 2) {$v_4$};
\node[draw,circle,red]		(D) at (2,0) {a};
\node[draw,circle,green]	(D) at (1,1) {b};
\node[draw,circle,blue] 	(D) at (2,1) {c};
\node[draw,circle,brown]	(D) at (0,3) {d};
\node[draw,circle,gray] 	(D) at (2,3) {e};
\end{tikzpicture}
\end{minipage}
	\begin{enumerate}
	\item[\textcolor{red}{a}]
	$v_5 = (v∖\{a_3\})∪\{x_1\}$.
	In this case $v_3,v_4,v_5$ all have distance $2$ from each other and thus by Lemma \ref{Diag} $A$ is not shattered.
	\item[\textcolor{green}{b}]
	$v_5 = (v∖\{a_2\})∪\{x_2\}$.
	Then $v_1,v_2,v_3,v_5$ form case $II$.
	\item[\textcolor{blue}{c}]
	$v_5 = (v∖\{a_3\})∪\{x_2\}$.
	In this case $v_1,v_4,v_5$ all have distance $2$ from each other and thus by Lemma \ref{Diag} $A$ is not shattered.
	\item[\textcolor{brown}{d}]
	$v_5 = (v∖\{a_1\})∪\{x_4\}$.
	In this case $v_1,v_3,v_4,v_5$ form case $VI$.
	\item[\textcolor{gray}{e}]
	$v_5 = (v∖\{a_3\})∪\{x_4\}$.
	In this case $v_3,v_4,v_5$ all have distance $2$ from each other and thus by Lemma \ref{Diag} $A$ is not shattered.
	\end{enumerate}

%% file: Johnson/Dimension/Cases/CaseC.tex
\begin{minipage}[t]{.5\textwidth}
\begin{align*}
	v_1 = (v∖\{a_1\})∪\{x_1\}\\
	v_2 = (v∖\{a_1\})∪\{x_2\}\\
	v_3 = (v∖\{a_1\})∪\{x_3\}\\
	v_4 = (v∖\{a_1\})∪\{x_4\}\\
\end{align*}
\end{minipage}
\begin{minipage}[t]{.5\textwidth}
\begin{tikzpicture}[yscale=-1,thick,baseline={([yshift=9ex]current bounding box.center)}]
\node[draw,circle] (A) at (0 , 0) {$v_1$};
\node[draw,circle] (B) at (0 , 1) {$v_2$};
\node[draw,circle] (C) at (0 , 2) {$v_3$};
\node[draw,circle] (D) at (0 , 3) {$v_4$};
\node[draw,circle,red]		(D) at (1,0) {a};
\node[draw,circle,green]	(D) at (0,4) {b};
\node[draw,circle,blue]		(D) at (1,4) {c};
\end{tikzpicture}
\end{minipage}
	\begin{enumerate}
	\item[\textcolor{red}{a}]
	$v_5 = (v∖\{a_2\})∪\{x_1\}$.
	Here $v_2,v_3,v_4,v_5$ form case $VI$. 
	\item[\textcolor{green}{b}]
	$v_5 = (v∖\{a_1\})∪\{x_5\}$.
	Let $w$ be such that $N(w)∩A=\{v_1,v_2,v_3\}$.
	Observe that $w≠v_4$ since $v_5∈N(v_4)$ so we will need an alternative $w$. 
	We have $2$ cases: either $d(v,w)=1$
	 or $d(v,w)=2$. 

	Let $w=(v∖\{a\})∪\{x\}$. Since we have to exclude $v_5$ from $N(w)$ then by Lemma \ref{d1} we cannot have $a=a_1$.
	So in order to have $v_1∈N(w)$ we must have $x=x_1$ but  in order to have $v_2∈N(w)$ we must have $x=x_2$, a contradiction.

	Let $w=(v∖\{a,b\})∪\{x,y\}$. 
	In order to have $v_1∈N(w),v_2∈N(w)$ and $v_3∈N(w)$ Lemma \ref{d2} gives us we must have
	$\{x_1,x_2,x_3\}⊆\{x,y\}$, a  contradiction. 	
	\item[\textcolor{blue}{c}]
	$v_5 = (v∖\{a_2\})∪\{x_5\}$	
	Here $v_2,v_3,v_4,v_5$ form case $VI$.
	\end{enumerate}

%% file: Johnson/Dimension/Cases/CaseE.tex
\begin{minipage}[t]{.5\textwidth}
\begin{align*}
	v_1 = (v∖\{a_1\})∪\{x_1\}\\
	v_2 = (v∖\{a_1\})∪\{x_2\}\\
	v_3 = (v∖\{a_2\})∪\{x_3\}\\
	v_4 = (v∖\{a_2\})∪\{x_4\}\\
\end{align*}
\end{minipage}
\begin{minipage}[t]{.5\textwidth}
\begin{tikzpicture}[yscale=-1,thick,baseline={([yshift=9ex]current bounding box.center)}]
\node[draw,circle] (A) at (0 , 0) {$v_1$};
\node[draw,circle] (B) at (0 , 1) {$v_2$};
\node[draw,circle] (C) at (1 , 2) {$v_3$};
\node[draw,circle] (D) at (1 , 3) {$v_4$};
\node[draw,circle,red]	(D) at (1,0) {a};
\node[draw,circle,blue] (D) at (2,0) {b};
\node[draw,circle,green](D) at (0,4) {c};
\node[draw,circle,brown] (D) at (2,4) {d};
\end{tikzpicture}
\end{minipage}
	\begin{enumerate}
	\item[\textcolor{red}{a}]
	$v_5 = (v∖\{a_2\})∪\{x_1\}$.
	Here $v_2,v_3,v_4,v_5$ form case $VI$ 
	\item[\textcolor{green}{b}]
	$v_5 = (v∖\{a_3\})∪\{x_1\}$.
	Here $v_2,v_3,v_5$ all have distance $2$ from each other and thus by Lemma \ref{Diag} $A$ is not shattered.
	\item[\textcolor{blue}{c}]
	$v_5 = (v∖\{a_1\})∪\{x_5\}$.
	In this case $v_1,v_2,v_3,v_5$ form case $VI$
	\item[\textcolor{brown}{d}]
	$v_5 = (v∖\{a_3\})∪\{x_5\}$.
	Here $v_1,v_3,v_5$ all have distance $2$ from each other and thus by Lemma \ref{Diag} $A$ is not shattered.
	\end{enumerate}

%% file: Johnson/Dimension/Tight.tex
\begin{theorem}

The \VC-dimension of the edge relation in the Johnson graph 
$J(m,k)$ is $4$ if and only if 
$1<k<m-1$ and 
$|V(J(m,k))|=\binom{m}{k}≥16$.

\begin{proof}
If $|V(J(m,k))|<16=2^4$ then the set system induced by the edge relation has fewer than $16$ sets.
Thus by the pigeonhole principle the \VC-dimension of the edge relation  is less than $4$.  

Assume $\binom{m}{k}≥16$ and $1<k<m-1$. 
Here we again  rely on $J(m-1,k-1)$ and $J(m-1,k)$ being induced subgraphs of $J(m,k)$. 
We also observe that $J(m,k)$ is isomorphic to $J(m,m-k)$.
So since $\binom{m}{k}≥16$ then either $J(7,2)$ or $J(6,3)$ are induced subgraphs of $J(m,k)$.

Since removing vertices from a graph can only decrease \VC-dimension it now suffices to show
that the edge relation has \VC-dimension $4$ in $J(7,2)$ and $J(6,3)$.
\ifConference
In Figure \ref{tightfig} in Appendix \ref{app:Tight} we show choices for vertices $v_1,v_2,v_3,v_4$ such that  $A=\{v_1,v_2,v_3,v_4\}$ is shattered by the edge relation, 
along with how each subset of $A$ can be obtained.
\else
In Figure \ref{tightfig}  we show choices for vertices $v_1,v_2,v_3,v_4$ such that  $A=\{v_1,v_2,v_3,v_4\}$ is shattered by the edge relation, 
along with how each subset of $A$ can be obtained.

\begin{figure}
\setlength{\columnseprule}{0.4pt}
\begin{multicols}{2}
\paragraph*{$J(7,2)$}
\begin{align*}
	v_1 =\{1,3\}\\
	v_2 =\{1,4\}\\
	v_3 =\{1,5\}\\
	v_4 =\{1,6\}\\
\end{align*}
\begin{align*}
A∩N(\{2,7\})	=&∅			\\
A∩N(\{3,7\})	=&\{v_1\}		\\
A∩N(\{4,7\})	=&\{v_2\}		\\
A∩N(\{5,7\})	=&\{v_3\}		\\
A∩N(\{6,7\})	=&\{v_4\}		\\
A∩N(\{3,4\})	=&\{v_1,v_2\}		\\
A∩N(\{3,5\})	=&\{v_1,v_3\}		\\
A∩N(\{3,6\})	=&\{v_1,v_4\}		\\
A∩N(\{4,5\})	=&\{v_2,v_3\}		\\
A∩N(\{4,6\})	=&\{v_2,v_4\}		\\
A∩N(\{5,6\})	=&\{v_3,v_4\}		\\
A∩N(v_4)	=&\{v_1,v_2,v_3\}	\\
A∩N(v_3)	=&\{v_1,v_2,v_4\}	\\
A∩N(v_2)	=&\{v_1,v_3,v_4\}	\\
A∩N(v_1)	=&\{v_2,v_3,v_4\}	\\
A∩N(\{1,2\})	=&A			\\
\end{align*}

\paragraph*{$J(6,3)$}
\begin{align*}
	v_1 =\{2,3,4\}\\
	v_2 =\{1,3,4\}\\
	v_3 =\{1,3,5\}\\
	v_4 =\{1,2,5\}\\
\end{align*}

\begin{align*}
A∩N(\{4,5,6\})	=&∅			\\
A∩N(\{2,3,6\})	=&\{v_1\}		\\
A∩N(\{v_1\})	=&\{v_2\}		\\
A∩N(\{v_4\})	=&\{v_3\}		\\
A∩N(\{1,2,6\})	=&\{v_4\}		\\
A∩N(\{3,4,6\})	=&\{v_1,v_2\}		\\
A∩N(v_2)	=&\{v_1,v_3\}		\\
A∩N(\{2,4,5\})	=&\{v_1,v_4\}		\\
A∩N(\{1,3,6\})	=&\{v_2,v_3\}		\\
A∩N(v_3)	=&\{v_2,v_4\}		\\
A∩N(\{1,5,6\})	=&\{v_3,v_4\}		\\
A∩N(\{3,4,5\})	=&\{v_1,v_2,v_3\}	\\
A∩N(\{1,2,4\})	=&\{v_1,v_2,v_4\}	\\
A∩N(\{2,3,5\})	=&\{v_1,v_3,v_4\}	\\
A∩N(\{1,4,5\})	=&\{v_2,v_3,v_4\}	\\
A∩N(\{1,2,3\})	=&A			\\
\end{align*}

\end{multicols}
\caption{Examples of shattered sets of size $4$ in $J(7,2)$ and $J(6,3)$}
\label{tightfig}
\end{figure}
\fi
So the $\VC$-dimension of the edge relation is at least $4$ in both $J(6,3)$ and $J(7,2)$.
This shows that the $\VC$-dimension of the edge relation is at least $4$ in all Johnson graphs $J(m,k)$ where $\binom{m}{k}≥16$ and $1<k<m-1$.
Theorem \ref{VC-dimTheorem} shows us that the edge relation has $\VC$-dimension at most $4$ in all Johnson graphs so this bound is tight whenever $\binom{m}{k}≥16$ and $1<k<m-1$.
\end{proof}

\end{theorem}

%% file: Johnson/Dimension/Corollary.tex
It is known that  boolean combinations of formulas with bounded \VC-dimension also have finite \VC-dimension (See Lemma 2.9 in \cite{SimonNIP}).
Since the only relations in the language of graphs are the edge relation and equality, and equality always has a \VC-dimension at most $1$ we get.
\begin{corollary}
Every quantifier free formula in the language of graphs has finite \VC-dimension on \Johnson.
\end{corollary}

%% file: Johnson/Density/NewArg.tex
\begin{theorem}
The \VC-density of the edge relation on $\Johnson$ is $2$.
\begin{proof}
First  we show that the \VC-density is at least $2$.
Assume without loss of generality that $m>2k$ and let $X$ be the underlying set of $J(m,k)$.
Fix a  vertex $v=\{a_i|1≤i≤k\}$ in $J(m,k)$, and let $(x_i)_{i=1}^{k}$ be distinct elements of $X$ such that $x_i∉v$ for all $i$.
Define  $A:=\{(v∖\{a_i\})∪\{x_i\}|1≤i≤k\}$.
Then for any pair of vertices $v_i:=(v∖\{a_i\})∪\{x_i\}$ and $v_j=(v∖\{a_j\})∪\{x_j\}$ we have that $N((v∖\{a_i\})∪\{x_j\})∩A=\{v_i,v_j\}$. 
There are 
$\frac{|A|^2-|A|}{2}$
such pairs so the $\VC$-density of the edge relation on $\Johnson$  is at least $2$.

Now we show that the $\VC$-density of the edge relation on $\Johnson$ is  at most $2$.
Let $A$  be a set of vertices in $J(m,k)$, and $π(n)$ be the shatter function for the edge relation on $J(m,k)$.
Let  $|A|=n$  and $A$ be maximally shattered by the edge relation for sets of size $n$.
Let
\begin{align*} 
S(A)&=\{N(u)∩A|u∈V(G)\},\\
C_1(A) &= \{ N∈ S(A): N \text{ is      a clique}\},\text{ and }\\
C_2(A) &= \{ N∈ S(A): N \text{ is not  a clique}\}.
\end{align*}
By our assumption that $A$ is maximally shattered we have  $|S(A)|= \pi(n)$.  
Note also that $S(A)=C_1(A)∪C_2(A)$ so we deal with those two cases separately.

$|C_1(A)|≤\frac{5|A|^2+3|A|}{2}$: 
There are at most $\frac{|A|^2+|A|}{2}$ cliques of size $2$ or less  in $S(A)$. 
There are at most $|A|$ cliques $C$ in $S(A)$ such that $C=A∩N(v)$ for some $v∈A$.

Now assume we have $C=A∩N(v)$ for some $v∉A$ and further assume that $|C|≥3$.
We want to show that then the clique $C$ is of the form $A∩Q$ for some maximal clique $Q$ of $J(m,k)$. 
We then argue that there can be at most $2|A|^2$ maximal cliques of $J(m,k)$ that intersect $A$ in more than one vertex.

Note that in any graph $G$ a  maximal clique $Q$ of $G$  is contained in $N(u)∪\{u\}$ for all $u∈Q$ so $Q∖\{u\}$ is a maximal clique in $G[N(u)]$.
It is easy to see that the maximal cliques of the rook's graph $R(m,k)$ are the rows and columns.
So by Lemma $\ref{d0}$ we find that for every vertex $u$ in $J(m,k)$ the maximal cliques of $J(m,k)$  that $u$ belongs to are of the form  $Z∪\{u\}$ where $Z$ is a  row or a  column of the rook's graph $J(m,k)[N(u)]$.

Since $|C|≥3$ we know by Lemma \ref{d1} the only vertices connected to all vertices in $C$  are $v$ and those vertices that share that row or column with all of $C$, in the rook's graph induced by $N(v)$, and therefore lie in $N(v)$.
It follows that  $C=A∩Q$ for some maximal clique $Q$ of $J(m,k)$.

For every vertex $u∈A$ we have that $A$ intersects at most $|A|$ rows and at most $|A|$ columns of the rook's graph induced by $N(u)$. 
So $u$ can be a member of at most $2|A|$ maximal cliques of $J(m,k)$ that intersect $A$ in more than two vertices.
So the number of maximal cliques of $J(m,k)$ that intersect $A$ in more than two vertices is at most $2|A|^2$.

$|C_2(A)|≤4|A|^2$: 
This holds since every pair of vertices at distance $2$ from each other can by Lemma \ref{JohnsonIntersect} be contained in the neighbourhood of at most $4$ vertices and there are at most $|A|^2$ such pairs.

So we get that 
$
|S(A)|	≤|C_1(A)|+|C_2(A)|
	≤\frac{5|A|^2+3|A|}{2}+4|A|^2
	=\frac{13|A|^2+3|A|}{2}∈\mathcal{O}(|A|^2)
$.
\end{proof}
\end{theorem}

%% file: Hamming/Main.tex
\section{Hamming Graphs}
\label{Hamming}
In this section we will give technical lemmas for dealing with Hamming graphs  and prove our main results  on the \VC-dimension and \VC-density of the edge relation in such graphs.

\input{Hamming/Lemmas/Main.tex}
\input{Hamming/Dimension.tex}

\input{Hamming/Density.tex}

%% file: Hamming/Lemmas/Main.tex
\input{Hamming/Lemmas/d0.tex}
\input{Hamming/Lemmas/d1.tex}
\input{Hamming/Lemmas/d2.tex}

Lemma \ref{Hd2} implies the following.
\input{Hamming/Lemmas/subdiv-clique.tex}

\input{Hamming/Lemmas/Neighborhoods.tex}


%% file: Hamming/Lemmas/d0.tex
\ifConference
\ifAppendix
\setcounter{theorem}{\value{Hd0}}
\else
\newcounter{Hd0}
\setcounter{Hd0}{\value{theorem}}
\fi
\fi
\begin{lemma}
\label{Hd0}
Let $v$ be a vertex in the Hamming graph $H(d,q)$.
Then $N(v)$ induces a disjoint union of $d$  copies of $K_{q-1}$.
\ifAppendix
\begin{proof}

We observe that for each coordinate $j$  the set of neighbours of $v$ that disagree with $v$ in the $j$-th coordinate  has size $q-1$ and since those vertices all agree in all but the $j$-th coordinate they form a clique.
If two vertices $u,w∈N(v)$ disagree with $v$ in different coordinates, say $i$  and $j$ respectively, then $u$ and $w$ disagree with each other in the $i$-th and the $j$-th coordinate and thus they are non-adjacent.
\end{proof}
\else
\label{Hd0}
\fi
\end{lemma}

%% file: Hamming/Lemmas/d1.tex
\begin{lemma}
\label{Hd1}
Let $u$ and $v$ be vertices in the Hamming Graph $H(d,q)$ with $d(u,v)=1$.
Let $1≤i≤d$ be such that $u$ and $v$ agree on all but the $i$-th coordinate. 
Then $N(u)∩N(v)$ is a clique of size $q-2$ whose members are all vertices  $w$ that agree with $u$ and $v$ in all but the $i$-th coordinate.
%
\ifAppendix
\begin{proof}

Since $u$ and $v$ are neighbours we know that they agree in all but one coordinate namely the $i$-th.
All vertices that agree with $u$ and $v$ on all coordinates except the $i$-th form a clique.
Since each coordinate can have $q$ different values there are $q-2$ such vertices that are neither $u$ nor $v$.

Take any vertex $w$ neighboring $u$ that agrees with $u$  on all coordinates except the $j$-th for a $j≠i$.
Now we know that $w$ and $u$ agree on the $i$-th coordinate but since $u$ and $v$ disagree on this coordinate we 
get that $w$ and $v$ disagree on the $i$-th and the $j$-th coordinate and therefore do not have an edge between them.
\end{proof}
\else
\label{Hd1}
\fi
\end{lemma}

%% file: Hamming/Lemmas/d2.tex
\begin{lemma}
\label{Hd2}
Let $u=(u_k)_{k=1}^d$ and $v=(v_k)_{k=1}^d$ be vertices in the Hamming Graph $H(d,q)$ with $d(u,v)=2$.
Let $1≤i<j≤d$ be such that $u_i≠v_i$, $u_j≠v_j$ and $u_k=v_k$ for every $k∉\{i,j\}$. 
Then $N(u)∩N(v)$ has exactly  two vertices and they are not connected, namely $x=(x_k)_{k=1}^d$ and $y=(y_k)_{k=1}^d$ where 
$x_i=u_i$ and $x_k=v_k$ for all $k≠i$,
and
$y_i=v_i$ and $y_k=u_k$ for all $k≠i$.
\ifAppendix
\begin{proof}
Since $u$ and $v$ disagree on both the $j$-th and the $i$-th coordinates any vertex $w∈N(u)∩N(v)$ 
will have to agree with $u$ on either the $i$-th or the $j$-th coordinate and with $v$ on the other one of those.

\end{proof}
\else
\label{Hd2}
\fi
\end{lemma}

%% file: Hamming/Lemmas/subdiv-clique.tex
\begin{corollary}
\label{cor:Hd2}
 The open 2-neighbourhood in the Hamming Graph $H(d,2)$ induces 
 the 1-subdivision of the complete graph $K_d$.
\end{corollary}
Since $H(d,2)$ is an induced subgraph of $H(d,q)$ for $q≥2$ it follows that $\Hamming$ has unbounded local clique-width as mentioned in the introduction.

%% file: Hamming/Lemmas/Neighborhoods.tex
\begin{lemma}
Let $u$ and $v$ be vertices in the Hamming graph $H(d,q)$ then

$$
|N(v)∩N(w)|=
	\begin{cases}
	d(q-1)	&\text{ if } d(u,v)=0\\
	q-2	&\text{ if } d(u,v)=1\\
	2	&\text{ if } d(u,v)=2\\
	0	&\text{ if } d(u,v)≥3
	\end{cases}	
$$
\begin{proof}
This follows immediately from Lemmas \ref{Hd0},\ref{Hd1}, and \ref{Hd2}.
\end{proof}
\end{lemma}

%% file: Hamming/Dimension.tex

\begin{theorem}
\label{HDim}
The \VC-dimension of the edge relation in a Hamming graph is at most $3$.
\begin{proof}
Assume there is a set $A'$ with $|A'|>3$ which is shattered by the edge relation.
Then there is a set $A⊆A'$ with $|A|=4$ which is shattered by the edge relation.
Let $A=\{v_1,v_2,v_3,v_4\}$,
let $v$ be such that $N(v)∩A=A$ and
 $w$ be a vertex such that $N(w)∩A=\{v_1,v_2,v_3\}$.
Since $v≠w$ and $|N(v)∩N(w)|>2$ we have that $d(v,w)=1$, so  the intersection of 
$N(v)$ and $N(w)$ is a clique.
Now we have two cases: either  $v_4=w$ or $v_4≠w$.

Assume 
$v_4=w$ so $A$ induces a clique.
Let $u$ be such that $N(u)∩A= \{v_1,v_2\}$.
Since $A$ is a clique we know that $u∉A$.
More importantly $u$ cannot belong to the copy of $K_{q-1}$ in  $N(v)$  that contains  $A$ so by  Lemma \ref{Hd0} $d(u,v)=2$.
But then $N(v)∩N(u)$ by Lemma \ref{Hd2} has two vertices that are not adjacent, in contradiction with $A$ being a clique.

Assume 
$v_4≠w$. 
Then we know that $v_4∉N(v)∩N(v_1)$ since otherwise it would be in $N(w)$ in contradiction with
$N(w)∩A=\{v_1,v_2,v_3\}$. 
Then $d(v_4,v_1)=2$ and  similarly  $d(v_4,v_2)=2$.
Let $u$ be a vertex such that $N(u)∩A=\{v_1,v_2,v_4\}$.
Since $u≠v$ and  $|N(u)∩N(v)|>2$ we have by Lemma  \ref{Hd1} that $N(u)∩A⊆N(u)∩N(v)$ is a clique,
in contradiction with  $d(v_4,v_1)=2$.

\end{proof}
\end{theorem}

\begin{theorem}
The \VC-dimension of the edge relation on the Hamming graph $H(d,q)$ is $3$ if and only if at least  one of the following holds.
\begin{enumerate} 
\item $d≥3$ and $q≥3$.
\item $d≥2$ and $q≥4$.
\item $d≥4$ and $q≥2$.
\end{enumerate}
\begin{proof}
Note that if $d≤d'$ and $q≤q'$ then $H(d,q)$ is an induced subgraph of $H(d,q)$.
Since removing vertices from a graph can only decrease \VC-dimension it now suffices to show
that the edge relation has \VC-dimension $3$ in $H(3,3)$,$H(2,4)$ and $H(4,2)$.

In Table \ref{HTightTable} we give examples of shattered sets $A=\{v_1,v_2,v_3\}$ in $H(2,4)$, $H(3,3)$ and $H(4,2)$.
The last $3$ columns in the second table show choices of $x$, in the different graphs, such that $A∩N(x)$ is the subset shown in the first column.

%
%
%
%
%
%

\begin{table}
\begin{center}
\begin{tabular}{c|c|c|c}
$A$	&$H(2,4)$&$H(3,3)$&$H(4,2)$\\
\hline
$v_1$&$(0,1)$&$(0,0,1)$&$(0,0,0,1)$\\
$v_2$&$(0,2)$&$(0,1,0)$&$(0,0,1,0)$\\
$v_3$&$(0,3)$&$(1,0,0)$&$(0,1,0,0)$
\end{tabular}
\vspace{5pt}
\begin{tabular}{ c |c |c |c }
$A∩N(x)$ 		& $H(2,4)$ 	& $H(3,3)$ 	& $H(4,2)$ 	\\
\hline 
$∅$  			& $(1,0)$	& $(1,1,1)$	& $(1,1,1,1)$	\\
$\{v_1\}$		& $(1,1)$	& $(0,0,2)$	& $(1,0,0,1)$	\\
$\{v_2\}$		& $(2,2)$	& $(0,2,0)$	& $(1,0,1,0)$	\\
$\{v_3\}$		& $(3,3)$	& $(2,0,0)$	& $(1,1,0,0)$	\\
$\{v_1,v_2\}$		& $(0,3)$	& $(0,1,1)$	& $(0,0,1,1)$	\\
$\{v_1,v_2\}$		& $(0,2)$	& $(0,0,1)$	& $(0,1,0,1)$	\\
$\{v_1,v_2\}$		& $(0,1)$	& $(1,1,0)$	& $(0,1,1,0)$	\\
$\{v_1,v_2,v_3\}$       & $(0,0)$	& $(0,0,0)$	& $(0,0,0,0)$	\\
\end{tabular}
\end{center}
\caption{Examples of shattered sets in $H(2,4)$, $H(3,3)$ and $H(4,2)$}
\label{HTightTable}
\end{table}

\end{proof}
\end{theorem}

%% file: Hamming/Density.tex

\begin{theorem}
The \VC-density of the edge relation on  $\Hamming$ is  $2$.

\begin{proof}

First we observe that for $d>1$ a set such that any two vertices agree on all but the first two  $2$ coordinates has the property that 
$∀u,v∈A∃w(A∩N(w)=\{u,v\})$ so the \VC-density is at least $2$.

We now  need to show that  $π_{\Hamming}(n)∈\mathcal{O}(n^2)$ where $π_{\Hamming}$ is the shatter function for the edge relation on \Hamming.
We do this by giving a bound on a recursive formula for $π(n)$ and showing that it has a $\mathcal{O}(n^2)$ closed form.

Let $A$ be a maximally shattered set of size $n$ in the Hamming graph $H(d,q)$.
Let $v∈A$ and $\mathcal{S}$ be the class of all neigbourhoods in $H(d,q)$. 
Let $S_1=\{A∩S|S∈\mathcal{S} ∧ v∈S\}$.
Let $S_2=\{A∩S|S∈\mathcal{S} ∧ v∉S\}$.
Note that $|S_1∪S_2|=π(n)$ and $|S_2|≤π(n-1)$.

Every member of $S_1$ is an intersection of $A$ with a  neighborhood of neighbour of $v$.
Let
\begin{align*}
D_0&=\{v\},&&
D_1=A∩N(v),\\
D_2&=\{u∈A|d(u,v)=2\}, \text{ and }&&
D_3=\{u∈A|d(u,v)>2\}
\end{align*}

Then $D_3$ intersects no member of $S_1$ by definition of $D_3$.
By Lemma \ref{Hd2} every element of $D_2$ can be a member of at most $2$ sets of $S_1$ thus the total number of distinct sets containing $v$ and intersecting $D_2$ is at most $2|D_2|<2n$.

Since we have counted all members of $S_1$ that intersect $D_2$, and no members of $S_1$ intersect $D_3$ we only have left to count those members of $S_1$  that are subsets of $D_0∪D_1$.
By Lemma \ref{Hd0},  $N(v)$ induces a disjoint union of $d$ copies of $K_{q-1}$. 
Let $(Q_i)_{i=1}^d$ be sets such that for each $i$,  $Q_i$ is the set of  all vertices $u∈D_1$ that disagree with  $v$ in the $i$-th coordinate.
Note that $D_1=⋃_{i=1}^dQ_i$ and any element of $S_1$ which is a subset of $D_0∪D_1$ is a subset of $D_0∪Q_i$ for some $i$.

Moreover every subset of $D_0∪Q_i$ that is an element of $S_1$ is either:
$(D_0∪Q_i)∖\{u\}$ for some $u∈Q_i$, 
or $D_0∪Q_i$, or $\{v\}$,
thus the number of distinct elements of $S_1$ contained in $D_0∪D_1$ is at most 
$$
∑_{i=1}^d |Q_i|+\min(d,n)+1=|D_1|+\min(d,n)+1≤n+\min(d,n)+1.$$

So we have
\begin{align*}
π(n)
&=|S_1|			+|S_2|
≤|S_1|			+π(n-1)
≤2|D_2|+n+\min(d,n)+1	+π(n-1)
\\&≤ 2n+n+n+1		+π(n-1)
≤ 4n+1+π(n-1).
\end{align*}

By induction we get that  $π(n)≤4n^2+n$ for all $n$ and thus $π(n)∈\mathcal{O}(n^2)$.
This tells us that the $\VC$-density is at most $2$.
We have thus demonstrated that the $\VC$-density of the edge relation on $\Hamming$ is at least $2$ and at most $2$ and conclude that  it must be $2$.
\end{proof}

\end{theorem}

%% file: Appendix.tex
\ifConference

\newpage
\appendix
\Appendixtrue
\begin{center}\huge\bf Appendix \end{center}
\section{Full definitions}\label{app:Def}
\input{Definitions.tex}

\section{Omitted proofs of lemmas}\label{app:Lemma}

\input{Lemmas.tex}
\section{Cases of Theorem \ref{VC-dimTheorem}}\label{app:Cases}
\input{Johnson/Dimension/Cases/Cases.tex}
\section{Examples of shattered sets}\label{app:Tight}

\input{Johnson/Dimension/Appendix.tex}

\input{Hamming/Appendix.tex}

\fi

%% file: Definitions.tex
Here we have the definitions used in this paper laid out for the benefit of the reader
\begin{definition}
A \emph{set system} is a pair $(X,\mathcal{S})$ consisting of a universe set $X$ and a family $\mathcal{S}$  of subsets of $X$. 
\end{definition}
Set systems are sometimes also referred to as hypergraphs or range spaces.

Note that for every partitioned first order formula $ϕ(x,y)$ in the language of graphs and  every graph $G$ we get a set system
$(V(G),\mathcal{S}_ϕ)$ where  $\mathcal{S}_ϕ=\{\{x∈(V(G)|ϕ(x,y)\}|y∈V(G)\}$.
In this paper we study the formula "$Exy$" i.e. one that states that $x$ and $y$ are adjacent. 
The set system obtained from this formula on a graph $G$ is $(G,\mathcal{S}_E)$ where $\mathcal{S}_E:=\{N(v)|v∈V(G)\})$.

\begin{definition}
Let $(X,\mathcal{S})$ be a set system and  $A\subseteq X$ be a set.
We say that $A$ is \emph{shattered} by $\mathcal{S}$ if the class of intersections of sets in $\mathcal{S}$ with $A$ is the full powerset of $A$ i.e.
$$
\forall a \subseteq A \ \exists S \in \mathcal{S}\ a = A\cap S
$$
\end{definition}

By abuse of notation we say a set $A$ of vertices in a graph $G$   is shattered by the edge relation if it is shattered in $(V(G),\mathcal{S}_E)$.

\begin{definition}
Let  $(X,\mathcal{S})$  be a set system. 
We define the \emph{shatter function} $\pi_{\mathcal{S}}: \setN \to \setN$ as: 
$$
\pi_\mathcal{S}(n) := \max  \{ | \{S\cap A : S\in \mathcal{S}\} | : A\subseteq X \wedge |A|= n\}
$$
\end{definition}

\begin{definition}
Let $(X,\mathcal{S})$  be a set system.
The \emph{\VC-Dimension} of $(X,\mathcal{S})$ is:
$$
\VC((X,\mathcal{S})) = \begin{cases}
\sup\{n \in \setN\cup\{\infty\}: X \text{ has a subset of size $n$ shattered by } \mathcal{S} \text{ }\}  &\text{ if } \mathcal{S}\neq\emptyset \\
-\infty &\text{ if } \mathcal{S}=\emptyset
\end{cases}
$$
\end{definition}
We observe that the shatter function is $2^n$ for $n$ at most the $\VC$-dimension of the set system but for any $n$ greater than the $\VC$-dimension it is bounded by a polynomial of degree bounded by the $\VC$-dimension.
This leads to the following definition.

\begin{definition}
Let $(X,\mathcal{S})$  be a set system. 
Then the \emph{\VC-density} of $(X,\mathcal{S})$ is:
$$
\vc(X,\mathcal{S})=\begin{cases} 
\inf\{r\in \setR^{+}:\pi_{\mathcal{S}}(n) \in \mathcal{O}(n^r)\} \text{ if } \VC(\mathcal{S}) < \infty \\
\infty \text{ otherwise}
\end{cases}
$$
\end{definition}

\subsection{Classes of finite set  systems}
This work focuses on classes of finite graphs, but the definitions of VC-dimension and VC-density assume a single structure.
In particular the definition of VC-density assumes the structure to be infinite.
In this section we give definitions that extend the definitions of VC-dimension, shatter function, and VC-density in a natural way to classes of finite models.

\begin{definition}
Let $\mathcal{C}$ be a class of finite set systems. 
The VC-dimension of $\mathcal{C}$  is then said to be:
$$VC(\mathcal{C})=\max \{VC(X,\mathcal{S}):(X,\mathcal{S})\in \mathcal{C}\}$$
\end{definition}

\begin{definition}

Let $\mathcal{C}$ be a class of finite set systems. 
The shatter function of $\mathcal{C}$ is then said to be: 
$$
\pi_{\mathcal{C}}(n) = \max\{\pi_{\mathcal{S}}(n)\ : (X,\mathcal{S})\in \mathcal{C}\}
$$
\end{definition}

\begin{definition}

Let $\mathcal{C} $  be a  class of finite set systems. Then the \emph{VC-density} of $\mathcal{C}$ is:
$$
vc(\mathcal{C})=\begin{cases} 
\inf\{r\in \setR^{+}:\pi_{\mathcal{C}}(n) \in O(n^r)\} \text{ if } VC(\mathcal{C}) < \infty \\
\infty \text{ otherwise}
\end{cases}
$$

\end{definition}

\input{Johnson/Introduction.tex}

\input{Hamming/Introduction.tex}

\input{Rook.tex}

%% file: Johnson/Introduction.tex

\begin{definition}
Let $n$ and $k$ be natural numbers $m≥k$.
The \emph{ Johnson graph} $J(m,k)$ is a graph whose vertices correspond to the $k$-element subsets of a set of size $m$ and two vertices are adjacent if their corresponding sets intersect in all but one element.
i.e. their symmetric difference has size $2$.
\end{definition}

\begin{lemma}[\cite{johnsondist}]

Let $u$ and $v$ be vertices in a Johnson graph.
Then the distance  $d(u,v)=\frac{|u\symdiff v|}{2}$
\TODO[cite]
\end{lemma}

%% file: Hamming/Introduction.tex
\begin{definition}
Let $d$ and $q$ be natural numbers and $S$ a set with $|S|=q$.
The \emph{ Hamming graph} $H(d,q)$ is a graph whose vertices correspond  to ordered elements of $S^d$
and 
two vertices are adjacent if they agree in all but one coordinate.
\end{definition}

%% file: Rook.tex
\begin{definition}
Let $m,n∈ℕ$.
The \emph{rook's graph} $R(m,n)$ is the cartesian product of $K_m$ and $K_n$.
We can represent every vertex by a tuple from $R\times C$ where $|R|=m$ and $|C|=n$
and two vertices $(i,j),(k,l)$ are adjacent if and only if $i=k$ or $j=l$.
\end{definition}

Note that $R(n,n)=H(2,n)$.

%% file: Lemmas.tex
\input{Johnson/Lemmas/d0.tex}
\input{Johnson/Lemmas/d1.tex}
\input{Johnson/Lemmas/d2.tex}
\input{Hamming/Lemmas/d0.tex}
\input{Hamming/Lemmas/d1.tex}
\input{Hamming/Lemmas/d2.tex}

%% file: Johnson/Dimension/Appendix.tex
\begin{figure}[h]
\setlength{\columnseprule}{0.4pt}
\begin{multicols}{2}
\paragraph{$J(7,2)$}
\begin{align*}
	v_1 =\{1,3\}\\
	v_2 =\{1,4\}\\
	v_3 =\{1,5\}\\
	v_4 =\{1,6\}\\
\end{align*}
\begin{align*}
A∩N(\{2,7\})	=&∅			\\
A∩N(\{3,7\})	=&\{v_1\}		\\
A∩N(\{4,7\})	=&\{v_2\}		\\
A∩N(\{5,7\})	=&\{v_3\}		\\
A∩N(\{6,7\})	=&\{v_4\}		\\
A∩N(\{3,4\})	=&\{v_1,v_2\}		\\
A∩N(\{3,5\})	=&\{v_1,v_3\}		\\
A∩N(\{3,6\})	=&\{v_1,v_4\}		\\
A∩N(\{4,5\})	=&\{v_2,v_3\}		\\
A∩N(\{4,6\})	=&\{v_2,v_4\}		\\
A∩N(\{5,6\})	=&\{v_3,v_4\}		\\
A∩N(v_4)	=&\{v_1,v_2,v_3\}	\\
A∩N(v_3)	=&\{v_1,v_2,v_4\}	\\
A∩N(v_2)	=&\{v_1,v_3,v_4\}	\\
A∩N(v_1)	=&\{v_2,v_3,v_4\}	\\
A∩N(\{1,2\})	=&A			\\
\end{align*}

\paragraph{$J(6,3)$}

\begin{align*}
	v_1 =\{2,3,4\}\\
	v_2 =\{1,3,4\}\\
	v_3 =\{1,3,5\}\\
	v_4 =\{1,2,5\}\\
\end{align*}

\begin{align*}
A∩N(\{4,5,6\})	=&∅			\\
A∩N(\{2,3,6\})	=&\{v_1\}		\\
A∩N(\{v_1\})	=&\{v_2\}		\\
A∩N(\{v_4\})	=&\{v_3\}		\\
A∩N(\{1,2,6\})	=&\{v_4\}		\\
A∩N(\{3,4,6\})	=&\{v_1,v_2\}		\\
A∩N(v_2)	=&\{v_1,v_3\}		\\
A∩N(\{2,4,5\})	=&\{v_1,v_4\}		\\
A∩N(\{1,3,6\})	=&\{v_2,v_3\}		\\
A∩N(v_3)	=&\{v_2,v_4\}		\\
A∩N(\{1,5,6\})	=&\{v_3,v_4\}		\\
A∩N(\{3,4,5\})	=&\{v_1,v_2,v_3\}	\\
A∩N(\{1,2,4\})	=&\{v_1,v_2,v_4\}	\\
A∩N(\{2,3,5\})	=&\{v_1,v_3,v_4\}	\\
A∩N(\{1,4,5\})	=&\{v_2,v_3,v_4\}	\\
A∩N(\{1,2,3\})	=&A			\\
\end{align*}

\end{multicols}
\caption{Examples of shattered sets of size $4$ in $J(7,2)$ and $J(6,3)$}
\label{tightfig}
\end{figure}

%% file: Hamming/Appendix.tex
\begin{figure}
\begin{multicols}{2}
\begin{align*}
	v_1 =(1,1,2)\\
	v_2 =(1,2,1)\\
	v_3 =(2,1,1)\\
\end{align*}

\begin{align*}
a∩n((2,2,2))	=&∅		\\
a∩n((1,1,3))	=&\{v_1\}	\\
a∩n((1,3,1))	=&\{v_2\}	\\
a∩n((3,1,1))	=&\{v_3\}	\\
a∩n((1,2,2))	=&\{v_1,v_2\}	\\
a∩n((2,1,2))	=&\{v_1,v_3\}	\\
a∩n((2,2,1))	=&\{v_2,v_3\}	\\
a∩n((1,1,1))	=&A		\\
\end{align*}

\end{multicols}
\caption{examples of shattered sets of size $3$ in $H(3,3)$}.
\label{Htightfig}
\end{figure}